\documentclass[12pt]{amsart}
\usepackage{amsfonts,amssymb,color}
\usepackage[mathscr]{eucal}
\usepackage{amsmath, amsthm}
\usepackage{mathrsfs}
\usepackage{amsbsy}
\usepackage{dsfont}
\usepackage{bbm}
\makeindex
\input xypic
\xyoption{all}
\begin{document}
\baselineskip = 5mm
\newcommand \lra {\longrightarrow}
\newcommand \hra {\hookrightarrow}
\newcommand \ZZ {{\mathbb Z}} 
\newcommand \NN {{\mathbb N}} 
\newcommand \QQ {{\mathbb Q}} 
\newcommand \RR {{\mathbb R}} 
\newcommand \CC {{\mathbb C}} 
\newcommand \bcA {{\mathscr A}}
\newcommand \bcB {{\mathscr B}}
\newcommand \bcC {{\mathscr C}}
\newcommand \bcD {{\mathscr D}}
\newcommand \bcE {{\mathscr E}}
\newcommand \bcF {{\mathscr F}}
\newcommand \bcI {{\mathscr I}}
\newcommand \bcJ {{\mathscr J}}
\newcommand \bcM {{\mathscr M}}
\newcommand \bcP {{\mathscr P}}
\newcommand \bcS {{\mathscr S}}
\newcommand \bcT {{\mathscr T}}
\newcommand \bcU {{\mathscr U}}
\newcommand \bcX {{\mathscr X}}
\newcommand \bcY {{\mathscr Y}}
\newcommand \bcZ {{\mathscr Z}}
\newcommand \C {{\mathscr C}}
\newcommand \im {{\rm im}}
\newcommand \Hom {{\rm Hom}}
\newcommand \colim {{{\rm colim}\, }} 
\newcommand \End {{\rm {End}}}
\newcommand \Aut {\rm Aut}
\newcommand \coker {{\rm {coker}}}
\newcommand \id {{\rm {id}}}
\newcommand \supp {{\rm {Supp}}\, }
\newcommand \CHM {{\bcC \! \bcM }}
\newcommand \SP {{\bcS \! \bcP }}
\newcommand \DM {{\mathscr D\! \mathscr M}}
\newcommand \MM {{\mathscr M\! \mathscr M}}
\newcommand \uno {{\mathbbm 1}}
\newcommand \Le {{\mathbbm L}}
\newcommand \ptr {{\pi _2^{\rm tr}}}
\newcommand \Ob {{\rm Ob}}
\newcommand \PR {{\mathbb P}} 
\newcommand \AF {{\mathbb A}} 
\newcommand \Spec {{\rm {Spec}}}
\newcommand \Pic {{\rm {Pic}}}
\newcommand \Alb {{\rm {Alb}}}
\newcommand \Jac {{\rm {Jac}}}
\newcommand \Sym {{\rm {Sym}}}
\newcommand \Corr {{CH}}
\newcommand \cha {{\rm {char}}}
\newcommand \tr {{\rm {tr}}}
\newcommand \trdeg {{\rm {tr.deg}}}
\newcommand \alg {{\rm {alg}}}
\newcommand \gm {{\mathfrak {m}}}
\newcommand \gp {{\mathfrak {p}}}
\def\blue {\color{blue}}
\def\red {\color{red}}
\def\green {\color{green}}
\newtheorem{theorem}{Theorem}
\newtheorem{lemma}[theorem]{Lemma}
\newtheorem{sublemma}[theorem]{Sublemma}
\newtheorem{corollary}[theorem]{Corollary}
\newtheorem{example}[theorem]{Example}
\newtheorem{exercise}[theorem]{Exersize}
\newtheorem{proposition}[theorem]{Proposition}
\newtheorem{remark}[theorem]{Remark}
\newtheorem{notation}[theorem]{Notation}
\newtheorem{definition}[theorem]{Definition}
\newtheorem{conjecture}[theorem]{Conjecture}
\newtheorem{claim}[theorem]{Claim}
\newenvironment{pf}{\par\noindent{\em Proof}.}{\hfill\framebox(6,6)
\par\medskip}
\title[Non-trivial elements in the Abel-Jacobi kernels]
{\bf Non-trivial elements in the Abel-Jacobi kernels of higher-dimensional varieties}
\author{S. Gorchinskiy, V. Guletski\u \i }

\date{14 August 2011}

\begin{abstract}
\noindent The purpose of this paper is to construct non-trivial elements in the Abel-Jacobi kernels in any codimension by specializing correspondences with non-trivial Hodge-theoretical invariants at points with different transcendence degrees over a subfield in $\CC $.
\end{abstract}





\maketitle

\section{Introduction}
\label{s-intro}

A standard approach in the study of algebraic cycles presumes the mapping of Chow groups into the Weil cohomology groups via so-called cycle-class maps. The kernels of such mappings can not be always easily understood, so that we send then homologically trivial cycle classes into the intermediate Jacobians via the Abel-Jacobi mappings. The kernels of the Abel-Jacobi maps are even harder to understand and, actually, they are quite mysterious. For example, if $X$ is a $K3$ surface, its Albanese kernel is unknown, in spite of the fact that such surfaces are intensively studied from many other points of view. If $X$ is a surface of general type with $p_g=0$ then Bloch's conjecture, which is a codimension $2$ part of the general Bloch-Beilinson philosophy, predicts that the kernel of the Abel-Jacobi map is zero, but we are too far from the understanding of this conjecture. Working over the algebraic closure $\bar \QQ $ of rational numbers we have the Bloch-Beilinson's conjecture saying that there are no non-trivial cycle classes in the Albanese kernel of any smooth projective variety over $\bar \QQ $, so that first non-trivial examples of elements in the Abel-Jacobi kernels are expected to be rational over an extension of $\QQ $ whose transcendence degree is at least one. The first example of such a cycle class was constructed in \cite{Schoen} by Schoen. He also refers to unpublished Nori's results along the same line. Schoen's algebraic cycle is an external product of two cycles on the product of two curves. The general case in dimension $2$ has been considered in \cite{GreenGriffithsParanjape} where the authors proved an existence of a non-trivial class in the Albanese kernel for an arbitrary surface over $\bar \QQ $ which is defined over a field of transcendence degree $1$. In higher dimension the problem was considered in \cite{RosenschonSaito} and \cite{Kerr}.

Our vision of non-trivial elements in the Abel-Jacobi kernel has been inspired by the papers of Bloch, Bloch-Srinivas, Green and Griffiths, and it can be illustrated on the level of surfaces. Let $X$ be a smooth projective surface over $\bar \QQ $, let $H^2_{\rm tr}(X_{\CC })$ be the transcendent part in the second cohomology of the complex surface $X_{\CC }$, and let $P$ and $Q$ be two closed points on $X_{\CC }$. Let also $t_P$ and $t_Q$ be the transcendence degrees of the residue fields, respectively, $\bar \QQ (P)$ and $\bar \QQ (Q)$ over $\bar \QQ $. Suppose that $t_P=2$ and $t_Q<2$. If $H^2_{\rm tr}(X_{\CC })$ is non-trivial, say when $X_{\CC }$ is a $K3$-surface, then the point $P$ cannot be rationally equivalent to the point $Q$. In contrast, if $H^2_{\rm tr}(X_{\CC })$ vanishes, then $Q$ can be rationally equivalent to $P$, as it happens for surfaces of general type with $p_g=0$ and with known Bloch's conjecture. In other words, $H^2_{\rm tr}(X_{\CC })$ plays the role of a ``potential" for zero-cycles with different transcendence degrees to be rationally equivalent to each other on surfaces, and the purpose of the present paper is to give an evidence for this insight in high dimension.

Respectively, the method we use to construct non-trivial cycle classes in the Abel-Jacobi kernel is not relevant to external product of cycles on two varieties, but rather rely upon the interplay between algebraic cycles defined over a field of non-zero transcendence degree and correspondences, see \cite{Bloch} and \cite{BlochSrinivas}. Correspondences act on Chow groups and cohomology, and as such they have more powerful cohomological and Hodge-theoretic invariants which allow detect non-triviality of the original cycle classes. At the same time, correspondences control the motivic information ``at large". Searching correspondences which act non-trivially on suitable pieces in the Hodge-decomposition, and specializing them at the generic point on one variable, we construct desirable non-trivial classes in the Abel-Jacobi kernels (Theorem \ref{non-triviality}). Notice that this principle appears in a hidden form in \cite{MorihikoSaito}, and without a proof in the paper by Green and Griffiths, \cite{GreenGriffiths}, which is the departing point in our considerations. Our main result (Theorem \ref{maintheorem}) says:

\medskip

\begin{itemize}

\item[]{}

{\it Let $X$ and $S$ be two irreducible smooth projective varieties over $k$, and let $\alpha $ be an algebraic cycle class modulo rational equivalence on the product $S\times X$. Assume $\alpha _{\CC }$ acts non-trivially on suitable Hodge pieces of the cohomology groups. Then the difference between specializations of $\alpha $ (modified by the Albanese projector) at the generic and a closed points on $S$ is always a non-trivial element in the Abel-Jacobi kernel of the variety $X_{\CC }$.}

\end{itemize}

\medskip

It means that, in any codimension, a substantial contribution in to the Abel-Jacobi kernels arises from algebraic cycles with different transcendence degrees, whose non-triviality depends upon the transcendental parts in Hodge structures. Respectively, the transcendence degrees of our cycles are always greater than zero. However, the Lefschetz hyperplane section theorem allows to make these transcendence degrees as low as possible, see Remark \ref{lower}. The Beilinson's conjecture says that there are no non-trivial algebraic cycles in the Abel-Jacobi kernels of transcendence degree zero.

The paper is organized as follows. First we recall the language of Chow motives and prove some easy motivic lemmas in Section \ref{factorization}. Section \ref{AJK} is devoted to recollections on the Abel-Jacobi kernels. Some important lemma (Lemma \ref{lemma-AJ}) is proved there. In Section \ref{mainresult} we prove Theorem \ref{non-triviality} and Theorem \ref{maintheorem}. In Section \ref{someexamples} we apply Theorem \ref{maintheorem} to the case of $K3$-surfaces, and then to the case of threefolds with known the Lefschetz conjecture, in order to demonstrate how the main result works in practice. Finally, in the last Section \ref{mainexample} we explain how to use Theorem \ref{non-triviality} in order to produce new non-trivial cycle classes in the Abel-Jacobi kernels using non-vanishing differential forms.

\bigskip

{\sc Acknowledgements.} The authors are thankful to F.Bogomolov for a helpful discussion about non-vanishing differential forms on algebraic varieties. The first author was partially supported by the grants RFBR 11-01-00145-a, NSh-4713.2010.1 and MK-4881.2011.1 and AG Laboratory GU-HSE, RF government grant, ag. 11 11.G34.31.0023.

\bigskip

\section{Motivic factorization}
\label{factorization}

Below we will intensively use the language of Chow motives. As this language is not yet absolutely standard, for the convenience of the reader, we first recall some basic definitions and fix notation on them. Then we will prove a few motivic lemmas and deduce one proposition on motivic factorization, which will be used in what follows.

Let $k$ be a field. For any algebraic scheme over $k$ let $CH^n(X)$ be the Chow group of codimension $n$ algebraic cycles on $X$ with coefficients in $\QQ $. Thus, $CH^n(X)$ is actually a $\QQ $-vector space. Chow groups with coefficients in $\ZZ $ will be denoted by $CH^n_{\ZZ }(X)$.

The category of Chow motives $\CHM $ over $k$ will be given in contravariant notation. That is, if $X$ and $Y$ are two smooth projective varieties over $k$, and $X=\cup _jX_j$ is the decomposition of $X$ into its connected components, then the group of correspondences of degree $m$ from $X$ to $Y$ is defined by the formula
  $$
  \Corr^m(X,Y)=\oplus _jCH^{e_j+m}(X_j\times Y)\; ,
  $$
where $e_j$ is the dimension of the component $X_j$. For example, given a regular morphism $f:X\to Y$ over $k$, the class of the transpose $\Gamma _f^t$ of its graph $\Gamma _f$ is in $\Corr ^0(X,Y)$. For any two correspondences $f\in \Corr ^n(X,Y)$ and $g\in \Corr ^m(Y,Z)$ their composition $g\circ f$ is defined by the standard formula
  $$
  g\circ f={p_{13}}_*(p_{12}^*(f)\cdot p_{23}^*(g))\; ,
  $$
in which the central dot denotes the intersection of cycle classes in Fulton's sense, \cite{Fulton}, and the projections are obvious. Notice that the composition \mbox{$g\circ f$} belongs to the group $\Corr ^{m+n}(X,Z)$, so that the group $\Corr ^0(X,X)$ is an associative algebra with respect to the above composition, and the class of the diagonal on $X\times X$ is the unit in this algebra. Objects in $\CHM $ are triples $(X,p,n)$ where $p$ is a projector on $X$, i.e. an idempotent in the algebra $\Corr ^0(X,X)$, and $n$ is an integer. For two motives $M=(X,p,m)$ and $N=(Y,q,n)$, the Hom-group of morphisms from $M$ to $N$ is defined by the formula
  $$
  \Hom (M,N)=q\circ\Corr^{n-m}(X,Y)\circ p\; .
  $$
Given a smooth projective variety $X$ over $k$ its motive $M(X)$ is defined by the diagonal $\Delta _X$ weighted by zero, $M(X)=(X,\Delta _X,0)$. For any morphism $f:X\to Y$ over $k$ the transpose $\Gamma _f^{t}$ represents a cycle class in $\Corr ^0(Y,X)$ denoted by $M(f)$. Then we have a contravariant functor $M:\SP \lra \CHM $ from the category $\SP $ of smooth projective varieties over $k$ to the category of Chow motives $\CHM $, also over $k$.

The category $\CHM $ is tensor,
  $$
  (X,p,m)\otimes (Y,q,n)=(X\times Y,p\otimes q,m+n)\; ,
  $$
the motive
  $$
  \uno =M(\Spec (k))=(\Spec (k),\Delta _{\Spec (k)},0)
  $$
is the unit of this tensor product, the triple
  $$
  \Le =(\Spec (k),\Delta _{\Spec (k)},-1)
  $$
is called the Lefschetz motive, and its geometrical meaning is provided by the isomorphism $M(\PR ^1)=\uno \oplus \Le $. For any integer $m$ we will write $\Le ^m$ for the $m$-fold tensor power $\Le ^{\otimes m}$. The duality in $\CHM $ is defined by the formula
  $$
  (X,p,m)^{\vee}=(X,p^{\rm t},d-m)\; ,
  $$
where $X$ is a smooth projective variety of pure dimension $d$ over $S$. The category $\CHM $ is rigid with respect to the above defined structures.

A degree $m$ correspondence $\beta $ from a smooth projective variety $X$ to a smooth projective variety $Y$ over $k$ acts on the corresponding Chow groups,
  $$
  \beta _*:CH^i(X)\lra CH^{i+m}(Y)\; ,
  $$
by the formula
  $$
  \beta _*(\alpha )=(p_2)_*(p_1^*(\alpha )\cdot \beta )\; ,
  $$
for any $\alpha \in CH^i(X)$.


\begin{lemma}
\label{lemma-surjectChow}
Let $f:X\to Y$ be a surjective proper morphism of varieties over $k$. Then for any $p\ge 0$, the push-forward homomorphism
  $$
  f_*:CH_p(X)\lra CH_p(Y)
  $$
is surjective.
\end{lemma}

\begin{pf}
Use the fact that any proper surjective morphism of algebraic varieties admits a multisection.
\end{pf}

\begin{lemma}
\label{factorization1}
Let $X$, $Y$ and $Z$ be three smooth projective varieties over $k$, $X$ equi-dimensional, and let $f:Y\to Z$ be a regular morphism over $k$. Let $\alpha $ be an element in $\Corr ^m(Y,X)$, i.e. a correspondence of degree $m$ from $Y$ to $X$. Consider the element $\gamma =(f\times \id _X)_*(\alpha )$ sitting in the group $\Corr ^m(Z,X)$. Then the following diagram commutes:
  $$
  \diagram
  M(Z) \ar[rr]^-{M(f)} \ar[ddrr]_-{\gamma } & & M(Y) \ar[dd]^-{\alpha } \\ \\
  & & M(X)\otimes \Le ^{-m}
  \enddiagram
  $$
\end{lemma}

\begin{pf}
Let $g:Y\to Z\times Y$ be a morphism defined by the morphisms $f$ and $\id _X$, so that $[\Gamma _f^t]=g_*[Y]$. Let $p_{12}$, $p_{23}$ and $p_{13}$ are projections of the product $Z\times Y\times X$ onto the corresponding factors. Then we see that
  $$
  p_{12}^*[\Gamma _f^t]=[\Gamma _f^t]\times X=g_*[Y]\times X=(g\times \id _X)_*[Y\times X]\; .
  $$
Therefore,
  $$
  \begin{array}{rcl}
  \alpha \circ [\Gamma _f^t]
  &=&
  {p_{13}}_*(p_{12}^*[\Gamma _f^t]\circ p_{23}^*(\alpha )) \\
  &=&
  {p_{13}}_*((g\times \id _X)_*[Y\times X]\cdot p_{23}^*(\alpha ))
  \end{array}
  $$
By the projection formula:
  $$
  \begin{array}{rcl}
  \alpha \circ [\Gamma _f^t]
  &=&
  {p_{13}}_*((g\times \id _X)_*[Y\times X]\cdot p_{23}^*(\alpha )) \\
  &=&
  {p_{13}}_*((g\times \id _X)_*((g\times \id _X)^*p_{23}^*(\alpha )\cdot [Y\times X])) \\
  &=&
  {p_{13}}_*((g\times \id _X)_*(g\times \id _X)^*p_{23}^*(\alpha )) \\
  &=&
  {p_{13}}_*((g\times \id _X)_*(p_{23}\circ (g\times \id _X))^*(\alpha ))
  \end{array}
  $$
But
  $$
  p_{23}\circ (g\times \id _X)=\id _{Y\times X}\; ,
  $$
so that
  $$
  \begin{array}{rcl}
  \alpha \circ [\Gamma _f^t]
  &=&
  {p_{13}}_*((g\times \id _X)_*(p_{23}\circ (g\times \id _X))^*(\alpha )) \\
  &=&
  {p_{13}}_*(g\times \id _X)_*(\alpha ) \\
  &=&
  (p_{13}\circ (g\times \id _X))_*(\alpha ) \\
  \end{array}
  $$
Since
  $$
  p_{13}\circ (g\times \id _X)=g\times \id _X\; ,
  $$
we obtain:
  $$
  \begin{array}{rcl}
  \alpha \circ [\Gamma _f^t]
  &=&
  (p_{13}\circ (g\times \id _X))_*(\alpha ) \\
  &=&
  (g\times \id _X)_*(\alpha )=\gamma \; .
  \end{array}
  $$
\end{pf}

\bigskip

Let $X$ and $Y$ be two equi-dimensional varieties over $k$. Following \cite{Barbieri Viale}, we will say that a correspondence $\alpha\in CH^m(X,Y)$ is {\it balanced on the left} (respectively, {\it on the right}) if there exists an equi-dimensional Zariski closed subscheme $Z\subset X$ with $\dim (Z)<\dim (X)$ and an algebraic cycle $\Gamma $ on $X\times Y$, such that $\alpha =[\Gamma ]$ in $CH^m(X,Y)$ and the support of $\Gamma $ is contained in $Z\times Y$ (respectively, in $X\times Z$). Such $Z$ will be called the {\it pans} of balancing. We say that a correspondence $\alpha\in CH^m(X,Y)$ is {\it balanced} if $\alpha=\alpha_1+\alpha_2$, where $\alpha_1$ is balanced on the left, and $\alpha_2$ is balanced on the right.

\begin{lemma}
\label{lemma-suppleft}
Let $X$ and $S$ be equidimensional smooth projective varieties over $k$, and let $\alpha \in CH^m(S,X)$. Then $\alpha $ is balanced on the left if and only if there exists an equidimensional smooth projective variety $Z$ over $k$ with $\dim (Z)<\dim (S)$, such that $\alpha $ factors through $M(Z)$, that is $\alpha$ is a composition
  $$
  M(S)\lra M(Z)\lra M(X)\otimes \Le ^{-m}\; .
  $$
\end{lemma}

\begin{pf}
Suppose $\alpha $ is balanced on the left by the pan $Z$. Let $s:\tilde Z\to Z$ be the resolution of singularities on $Z$. Since $s\times \id _X$ is a surjective morphism, the induced homomorphism
  $$
  (s\times \id _X)_*:CH_d(\tilde Z\times X)\lra CH_d(Z\times X)
  $$
is surjective too, where $d$ is the dimension of $X$. Let $\tilde \gamma $ be a cycle class in $CH_d(\tilde Z\times X)$, such that $(s\times \id _X)_*\tilde \gamma =\gamma $. It exists by Lemma \ref{lemma-surjectChow}. Let $f=i\circ s$. Then
  $$
  (f\times \id _X)_*\tilde \gamma =
  (i\times \id _X)_*(s\times \id _X)_*\tilde \gamma =
  (i\times \id _X)_*\gamma =\pi .
  $$
By Lemma \ref{factorization1}, the following diagram is commutative:
  $$
  \diagram
  M(S) \ar[rr]^-{M(f)} \ar[ddrr]_-{\alpha } & & M(\tilde Z)
  \ar[dd]^-{\tilde \gamma } \\ \\
  & & M(X)\otimes \Le ^{-m}
  \enddiagram
  $$

\bigskip

Let us prove the converse. Suppose $\alpha$ factors through the Chow motive $M(Z)$,
  $$
  M(S)\stackrel{\beta }{\lra }M(Z)\stackrel{\gamma }{\lra }M(X)\otimes \Le ^{-m}\; ,
  $$
with $Z$ as in the statement of the lemma. Let
  $$
  W=\sum_i n_iW_i
  $$
be a cycle on $S\times Z$, such that $[W]=\beta $. Then $\dim (W_i)=\dim (Z)$ and, therefore, $\dim (p_S(W_i))<\dim (S)$, where $p_S:S\times Z\to S$ is the natural projection. This shows that there exists a subvariety $T\subset S$ with $\dim (T)<\dim (S)$, such that the support of $W$ is contained in $T\times Z\subset S\times Z$. Let $f:\widetilde T\to S$ be the composition of a resolution of singularities $\widetilde T\to T$ and the closed embedding $T\to S$.

Notice that the element
  $$
  \beta=[W]\in CH^0(S,Z)=CH_{\dim (Z)}(S\times Z)
  $$
is the push-forward of an element in $CH_{\dim (Z)}(T\times Z)$ with respect to the closed embedding $T\times Z\to S\times Z$. In addition, the natural morphism $\widetilde T\times Z\to T\times Z$ is surjective and proper. Hence, by Lemma~\ref{lemma-surjectChow}, there exists an element
$$
\tilde\beta\in CH^0(\widetilde T,Z)=CH_{\dim (Z)}(\widetilde T\times Z)\; ,
$$
such that $(f\times \id _Z)_*(\tilde \beta )=\beta $. By Lemma \ref{factorization1}, we have that
  $$
  \beta =\tilde \beta \circ M(f)\; .
  $$
Consequently,
  $$
  \alpha =\gamma \circ \beta =(\gamma \circ \tilde \beta )\circ M(f)\; .
  $$
Again by Lemma \ref{factorization1}, one has
  $$
  (\gamma \circ \tilde \beta )\circ M(f)=
  (f\times \id _X)_*(\gamma \circ \tilde \beta )\; ,
  $$
where $\gamma \circ \tilde \beta \in CH^m(\widetilde T,X)$. Since all elements in the image of the push-forward map
  $$
  (f\times \id _X)_*:CH^m(\widetilde T,X)\lra CH^m(S,X)
  $$
have representatives with support in $T\times X$, we obtain that $\alpha =(f\times \id _X)_*(\gamma \circ \tilde \beta )$ is balanced on the left.
\end{pf}

\begin{corollary}
\label{corol-suppright}
Let $X$ and $S$ be equidimensional smooth projective varieties over $k$, and let $\alpha \in CH^m(S,X)$. Then $\alpha $ is balanced on the right if and only if there exists an equidimensional smooth projective variety $Z$ over $k$ with $n:=\dim(X)-\dim(Z)>0$ such that $\alpha $ factors through $M(Z)\otimes \Le^{n-m}$, that is, $\alpha $ is a composition
$$
M(S)\lra M(Z)\otimes \Le ^{n-m}\lra M(X)\otimes \Le^{-m}\; .
$$
\end{corollary}

\begin{pf}
The corollary follows from Lemma~\ref{lemma-suppleft} by taking transpositions of the involved correspondences.
\end{pf}

\begin{remark}
\label{balancing}
{\rm
Certainly, one can now easily deduce a general result saying that $\alpha $ is balanced if and only if there exists equidimensional smooth projective varieties $Z$ and $T$ over $k$, such that $\alpha $ factors,
$$
M(S)\lra M(T)\oplus (M(Z)\otimes \Le ^{n-m})\lra M(X)\otimes \Le^{-m}\; ,
$$
and the difference $n=\dim(X)-\dim(Z)$ is greater than $0$. But we will not use it in the paper.
}
\end{remark}

\section{Abel - Jacobi kernels}
\label{AJK}

Here we recall some basic facts about the Abel-Jacobi maps and then prove some lemma about specialization of correspondences lying in the Abel-Jacobi kernels.

Let $X$ be an irreducible smooth projective variety over $\CC $, and let $H^i(X,\ZZ )$ be the $i$-th Betti cohomology group with coefficients in $\ZZ $. For any $p\geq 0$ let $CH^p_{\ZZ}(X)_0$ be the kernel of the cycle class homomorphism
  $$
  cl : CH^p_{\ZZ }(X)\lra H^{2p}(X,\ZZ )\; .
  $$
Extending coefficients from $\ZZ $ to $\QQ $ we obtain a rational cycle class map
  $$
  cl : CH^p(X)\lra H^{2p}(X,\QQ)\; ,
  $$
whose kernel will be denoted by $CH^p(X)_0$.

The Betti cohomology groups $H^i(X,\QQ )$ carry a pure Hodge structure. It means that there is a decreasing Hodge filtration $F^p$ on the $\CC $-vector space
  $$
  H^i(X,\CC )\cong H^i(X,\QQ )\otimes _{\QQ }\CC \; ,
  $$
compatible in a sense with the complex conjugate filtration $\bar F^p$. Let $H^{p,q}(X)$ be the adjoint quotient $(F^{p}/F^{p+1})H^{p+q}(X,\CC )$.

Consider two irreducible smooth projective complex algebraic varieties $X$ and $S$ and a correspondence $\alpha \in CH^m(S,X)$. The correspondence $\alpha $ defines a map
  $$
  \alpha _*:H^i(S,\QQ )\lra H^{i+2m}(X,\QQ )\; ,
  $$
which shifts the Hodge filtration by the formula
  $$
  \alpha _*(F^pH^i(S,\CC ))\subset F^{p+m}H^{i+2m}(X,\CC )\; .
  $$

\medskip

Let us describe this in some more detail. Let
  $$
  p_X : S\times X\to X\quad \hbox{and}\quad p_S:S\times X\to S
  $$
denote the natural projections. For any element $h\in F^pH^i(S,\CC)$, we have that
  $$
  p_S^*(h)\in F^pH^i(S\times X,\CC )\; ,
  $$
because pull-backs of elements in cohomology groups preserve the Hodge filtration, see \cite[Section 7.3.2]{Voisin}. Let $d$ be the dimension of the variety $S$. Since $\alpha $ is a correspondence of degree $m$ from $S$ to $X$, we obtain that
  $$
  cl(\alpha )\in F^{d+m}H^{2d+2m}(S\times X,\CC )\; ,
  $$
see Proposition 11.20 in loc. cit. Therefore,
  $$
  cl(\alpha )\cdot p_S^*(h)\in F^{d+m+p}H^{2d+2m+i}(S\times X,\CC )\; .
  $$
It follows that
  $$
  (p_X)_*(cl(\alpha )\cdot p_S^*(h))\in F^{m+p}H^{2m+i}(X,\CC )\; ,
  $$
as push-forwards of elements in cohomology groups shifts the Hodge filtration by the dimension of fibres, see again Section 7.3.2 in \cite{Voisin}. Thus, the action $\alpha _*$ induces maps also on adjoint quotients:
  $$
  \alpha _*:H^{p,q}(S)\lra H^{p+m,q+m}(X)\; .
  $$

\medskip

A $p$-th intermediate Jacobian is a compact complex torus
defined by the formula
$$
  J^{2p-1}(X)=H^{2p-1}(X,\CC )/({\rm Im}(H^{2p-1}(X,\ZZ ))+F^pH^{2p-1}(X,\CC ))\; ,
$$
where ${\rm Im}(H^{2p-1}(X,\ZZ))$ is the image of the natural map
  $$
  H^{2p-1}(X,\ZZ )\lra H^{2p-1}(X,\CC )\; ,
  $$
see \cite[Section 12.1.1]{Voisin}.
There is a special group homomorphism
  $$
  AJ : CH^p_{\ZZ }(X)_0\lra J^{2p-1}(X)\; ,
  $$
called the Abel-Jacobi map, whose precise description can be found in Section 12.1.2 of \cite{Voisin}. We will also work with the group
  $$
  J^{2p-1}(X)_{\QQ }=J^{2p-1}(X)\otimes_{\ZZ }\QQ \;
  $$
which can be also described as a quotient
  $$
  J^{2p-1}(X)_\QQ=H^{2p-1}(X,\CC )/((H^{2p-1}(X,\QQ)+F^pH^{2p-1}(X,\CC ))\; .
  $$
Extending coefficients of the integral Abel-Jacobi map from $\ZZ $ to $\QQ $, we obtain the rational Abel-Jacobi map
  $$
  AJ : CH^p(X)_0\lra J^{2p-1}(X)_{\QQ }\; ,
  $$
whose kernel we denote by $T^p(X)$. Thus, an element $\alpha \in CH^p_{\ZZ }(X)$ is in $T^{p}(X)$ if and only if $cl(n\cdot \alpha )=0$ and $AJ(n\cdot \alpha )=0$ for an integer $n$.

Notice that the action of correspondences on Chow groups is compatible with their actions on Hodge structures via the cycle class map and the Abel-Jacobi map.

\medskip

Denote by $J^{2p-1}(X)_{\alg }$ the largest complex subtorus of $J^{2p-1}(X)$ whose tangent space is contained in $H^{p-1,p}(X)$, see \cite[Section 12.2.2]{Voisin}. It follows $J^{2p-1}(X)_{\alg}$ is an abelian variety over $\CC$, loc.cit. Denote by $CH^p(X)_{\alg }$ (respectively, $CH^p_{\ZZ }(X)_{\alg }$) the group of algebraically trivial cycles with rational (respectively, integral) coefficients modulo rational equivalence. Then $CH^p(X)_{\alg }\subset CH^p(X)_0$, i.e.
  $$
  cl(CH^p(X)_{\alg })=0\; .
  $$
It is not hard to show that
  $$
  AJ(CH^p_{\ZZ }(X)_{\alg })\subset J^{2p-1}(X)_{\alg }
  $$
and
  $$
  AJ(CH^p(X)_{\alg })\subset (J^{2p-1}(X)_{\alg })_{\QQ}\; ,
  $$
loc.cit.

\bigskip

{\it From now on we assume that $k$ is an algebraically closed subfield in $\CC $.}

\medskip

Let $X$ be an irreducible smooth projective variety over $k$, and let $L$ be a finitely generated field extension of $k$. Let $\alpha $ be an element in the Chow group $CH^p(X_L)$, where $X_L=X\times _{\Spec (k)}\Spec (L)$. For any prime $l$ one has the \'etale cycle class homomorphism
  $$
  cl_l : CH^p(X_L)\lra
  H^{2p}_{\acute e t}(X_{\bar L},\QQ_l(p))\; .
  $$

Notice that the condition $cl_l(\alpha )=0$ can be also determined via Betti cohomology by choosing an embedding $\sigma :L\hra \CC $ over $k$. Indeed, any such $\sigma $ gives rise to an embedding $\bar \sigma :\bar L\hra \CC $ of an algebraic closure $\bar L$ of $L$ into $\CC $ over $\sigma $. Then we have a commutative diagram
  $$
  \diagram
  CH^p(X_L)\ar[dd]_-{\sigma ^*} \ar[rr]^-{cl_l} & & H^{2p}_{\acute e t}(X_{\bar L},\QQ_l(p)) \ar[dd]^-{\bar \sigma ^*} \\ \\
  CH^p(X_{\CC })\ar[rr]^-{cl_l} & & H^{2p}_{\acute e t}(X_{\CC },\QQ_l(p))
  \enddiagram
  $$
The right vertical homomorphism is an isomorphism, and $H^{2p}_{\acute e t}(X_{\CC },\QQ_l(p))$ is isomorphic to the Betti cohomology group $H^{2p}(X_{\CC },\QQ _l)$ by the well known properties of \'etale cohomology groups. Therefore, the condition $cl_l(\alpha )=0$ is equivalent to the condition $cl(\sigma ^*(\alpha ))=0$. Hence, the latter condition does not depend on the choice of the embedding $\sigma $, so that we can write
  $$
  cl(\alpha ) = 0 \; .
  $$
Let
  $$
  CH^p(X_L)_0 = \{ \alpha \in CH^p(X_L)\; |\; cl(\alpha )=0 \}
  $$
be a subgroup generated by homologically trivial cycle classes in $CH^p(X_L)$.

Let $\alpha $ be a cycle class in $CH^p(X_L)_0$. Then $\sigma ^*(\alpha )$ is in the kernel $CH^p(X_{\CC })_0$ of the complex cycle class map $cl:CH^p(X_{\CC })\to H^{2p}(X_{\CC },\QQ )$, and we can apply the Abel-Jacobi map to $\sigma ^*(\alpha )$. We will say that Abel-Jacobi class of $\alpha $ is zero, and write
  $$
  AJ(\alpha )=0\; ,
  $$
if $AJ(\sigma ^*(\alpha ))=0$ for any embedding $\sigma $ of $L$ into $\CC $ over $k$. Let
  $$
  T^p(X_L)=\{ \alpha \in CH^p(X_L)\; | \; AJ(\alpha )=0\} \; .
  $$
The group $T^p(X_L)$ is nothing but the Abel-Jacobi kernel for $X_L$.

\bigskip

Let $X$ and $S$ be two irreducible varieties over $k$, and assume that $X$ is smooth and projective. For any schematic point $x\in S$ let $k(x)$ be the residue field of $x$, as a point of a scheme over $k$. For any cycle class
  $$
  \alpha \in CH^p(S\times X)
  $$
let
  $$
  \alpha _x\in CH^p(X_{k(x)})
  $$
be the image of $\alpha $ under the pulling-back homomorphism
  $$
  CH^p(S\times X)\lra CH^p(X_{k(x)})
  $$
induced by the morphism
  $$
  \Spec (k(x))\times X \lra S\times X\; .
  $$
In particular, the generic point
  $$
  \eta =\Spec (k(S))
  $$
of the scheme $S$ gives rise to a restriction
  $$
  \alpha _{\eta }\in CH^p(X_{k(S)})
  $$
of the cycle class $\alpha $ at $\eta $ from the left.

\medskip

Let again $\alpha $ be a cycle class in the group $CH^p(S\times X)$. By
  $$
  \alpha _{\CC }\in CH^p((S\times X)_{\CC })
  $$
we denote the extension of scalars for the class $\alpha $ from $k$ to $\CC $. That is, $\alpha _{\CC }$ is the pull-back of $\alpha $ with respect to the natural flat morphism of schemes
  $$
  (S\times X)_{\CC }\to S\times X\; .
  $$
Since the scheme $(S\times X)_{\CC }$ is the same as the scheme $S_{\CC }\times X_{\CC }$, the class $\alpha _{\CC }$ is a correspondence from $S_{\CC }$ to $X_{\CC }$.

\bigskip

\begin{lemma}
\label{lemma-AJ}
For any element $\alpha \in CH^p(S\times X)$, we have that:

\medskip

  \begin{enumerate}

  \item[]{(i)}
  $cl(\alpha _{\eta })=0$ if and only if $cl((\alpha _{\CC })_x)=0$ for any closed point $x\in S_{\CC }$;

  \medskip

  \item[]{(ii)}
  assuming one of the equivalent conditions in $(i)$, $AJ(\alpha _{\eta })=0$ if and only if $AJ((\alpha _{\CC } )_x)=0$ for any closed point $x$ in $S_{\CC }$.

  \end{enumerate}

\medskip

\end{lemma}

\begin{pf}
To prove (i) recall that $cl(\alpha _{\eta })=0$ is equivalent to say that $cl((\alpha _{\eta })_{\CC })=0$ for a certain embedding of the function field $k(S)$ into $\CC $ over $k$. Suppose $cl((\alpha _{\CC })_x)=0$ for any closed point $x\in S_{\CC }$. Any embedding $k(S)\hra \CC $ over $k$ gives rise to a closed point $\eta _{\CC }\in S_{\CC }$, so that, in particular, $cl((\alpha _{\CC })_{\eta _{\CC }})=0$. The extension of scalars satisfies
  $$
  (\alpha _{\CC })_{\eta _{\CC }} = (\alpha _\eta )_{\CC }\; .
  $$
Therefore,
  $$
  cl((\alpha _{\eta })_{\CC })=cl((\alpha _{\CC })_{\eta _{\CC }})=0\; .
  $$

Now suppose $cl((\alpha _{\eta })_{\CC })=0$ under certain embedding of $k(S)$ into $\CC $ over $k$. Since $(\alpha _{\eta })_{\CC }$ is the same as $(\alpha _{\CC })_{\eta _{\CC }}$, we have that $cl((\alpha _{\CC })_{\eta _{\CC }})=0$. For any other closed point $x$ in $S_{\CC }$, different from $\eta _{\CC }$, the two specializations $(\alpha _{\CC })_x$ and $(\alpha _{\CC })_{\eta _{\CC}}$ are algebraically, and so homologically, equivalent to each other. Therefore,
  $$
  cl((\alpha _{\CC })_x)=cl((\alpha _{\CC })_{\eta _{\CC }})=0\; ,
  $$
which finishes the proof of $(i)$.

\medskip

Now let us prove $(ii)$. Assume $AJ((\alpha _{\CC })_x)=0$ for any closed point $x$ in $S_{\CC }$. In particular, $AJ((\alpha _{\CC })_{\eta _{\CC }})=0$, where $\eta _{\CC }$ is the closed point on $S_{\CC }$ induced by an embedding of $k(S)$ in to $\CC $ over $k$. As above, $(\alpha _\eta )_{\CC }=(\alpha _{\CC })_{\eta _{\CC }}$, so that $AJ((\alpha _{\eta })_{\CC })=0$. As it happens with respect to an arbitrary embedding of $k(S)$ in to $\CC $ over $k$, it means that $AJ(\alpha _{\eta })=0$.

Thus, it remains only to show the implication ``from the left to the right" in (ii). So, suppose $AJ(\alpha _{\eta })=0$, i.e. $AJ(\sigma ^*(\alpha _{\eta }))=0$ for any embedding $\sigma $ of the function field $k(S)$ into $\CC $ over $k$.

Let $n$ be the least common multiple of the denominators in the coefficients of a certain representative of $\alpha $. Then $\beta =n\cdot \alpha $ is in $CH^p_{\ZZ }(S\times X)$, and $AJ(\sigma ^*(\beta _{\eta }))$ is torsion in for any embedding $\sigma $ of $k(S)$ into $\CC $.

As the action of correspondences on Chow groups is compatible with their action on cohomology and their Hodge pieces, one has a commutative square
  $$
  \diagram
  CH^d_{\ZZ }(S_{\CC })_0 \ar[dd]_-{AJ}
  \ar[rr]^-{(\beta _{\CC })_*} & &
  CH^p_{\ZZ }(X_{\CC })_0 \ar[dd]^-{AJ} \\ \\
  \Alb (S_{\CC }) \ar[rr]^-{(\beta _{\CC })_*} & & J^{2p-1}(X_{\CC })_{\alg }
  \enddiagram
  $$
The lower horizontal $(\beta _{\CC })_*$ is a regular map between abelian varieties over $\CC $. Choosing an embedding $k(S)\hookrightarrow \CC$, let $\eta_{\,\CC}$ be the corresponding closed point on $S_{\CC}$, and let
  $$
  a : S_{\CC}\lra \Alb(S_{\CC})
  $$
be the Albanese morphism
  $$
  x\mapsto \int _{\eta_{\,\CC}}^x
  $$
induced by integration from $\eta_{\,\CC}$ to $x$ on $S_{\CC }$. As the above diagram commutes, a composition
  $$
  f=(\beta _{\CC })_*\circ a : S_{\CC}\lra J^{2p-1}(X_{\CC })_{\alg }
  $$
is a regular map sending any closed point $x\in S_{\CC }$ to the Abel-Jacobi class of the difference $(\beta _{\CC })_x-(\beta _{\CC })_{\eta _{\CC }}$. As we have shown above, the element $AJ((\beta _{\CC })_{\eta _{\CC }})$ is torsion, whence
  $$
  m\cdot f(x)=m\cdot AJ((\beta _{\CC })_x)
  $$
for some non-trivial integer $m$.

As $f$ is a regular map between two algebraic varieties, there exists a finitely generated field subextension
  $$
  k \subset K \subset \CC \; ,
  $$
an abelian variety $A$ over $K$, and a morphism of varieties over $K$
  $$
  g:S_K\lra A\; ,
  $$
such that $A_{\CC }=J^{2p-1}(X)_{\alg }$ and $g_{\CC }=f$.

Let
  $$
  \xi : \Spec(K(S_K))\lra S_K
  $$
be the generic point of $S_K$ as a scheme over $K$. This $\xi $  can be also obtained from the generic point $\eta $ of $S$ over $k$ by extending scalars from $k$ to $K$. Choose an embedding of $K(S_K)$ into $\CC $ over $K$. Automatically, it fixes an embedding
  $$
  k(S)\hra \CC \; ,
  $$
which may be different from the previous one, so that the corresponding closed point $\eta _{\CC }$ is now different from the previous $\eta _{\CC }$.

The new scalar extensions satisfy
  $$
  \xi_{\CC }=\eta _{\CC }\; .
  $$
As $\eta _{\CC }$ is a closed point on $S_{\CC }$,
  $$
  m\cdot f(\eta _{\CC })=
  m\cdot AJ((\beta _{\CC })_{\eta _{\CC }})\; ,
  $$
and the extension of scalars also gives
  $$
  AJ((\beta _{\CC })_{\eta _{\CC }})=
  AJ((\beta _\eta )_{\CC })\; .
  $$
Therefore,
  $$
  m\cdot g_{\CC }(\xi _{\CC })=m\cdot f(\xi _{\CC })=
  m\cdot f(\eta _{\CC })=
  m\cdot AJ((\beta _{\eta })_{\CC })\; .
  $$
Since $AJ(\alpha _{\eta })=0$, we have that $AJ((\beta _\eta)_{\CC })$ is torsion, with respect to the above embedding of the function fields in to $\CC $. Let $m'$ be a non-trivial integer, such that
  $$
  m'\cdot AJ((\beta _\eta)_{\CC })=0\; .
  $$

Then, $mm'\cdot g_{\CC }(\xi _{\CC })=0$ and so
  $$
  mm'\cdot g(\xi )=0\; ,
  $$
i.e. the morphism of varieties $g:S_K\to A$ over $K$
sends the generic point $\xi $ of $S_K$ to a schematic point lying in the $mm'$-torsion subgroup $A_{mm'}$ of $A$. Therefore, $g$ sends the whole variety $S_K$ into $A_{mm'}$. As $S_K$ is irreducible,
  $$
  g(S_K)=0\; .
  $$

It gives that
  $$
  f(S_{\CC })=0\; ,
  $$
because $f=g_{\CC }$. This means that
  $$
  AJ((\beta \,_{\CC })_x)=0
  $$
for any closed point $x$ on $S_{\CC }$. Since $\beta =n\cdot \alpha $, we see that
  $$
  AJ((n\cdot \alpha \,_{\CC })_x)=0
  $$
for any point $x\in S_{\CC }$. Therefore,
  $$
  n\cdot AJ((\alpha\,_{\CC })_x)=0
  $$
for any point $x\in S_{\CC }$.

Since the rational Abel-Jacobi map
  $$
  AJ : CH^p(X_{\CC })\lra J^{2p-1}(X_{\CC })_{\QQ }
  $$
takes its values in the $\QQ $-vector space $J^{2p-1}(X_{\CC })_{\QQ }$, we obtain that
  $$
  AJ((\alpha \,_{\CC })_x)=0
  $$
for any point $x\in S_{\CC }$.
\end{pf}

\begin{remark}
{\rm A similar argument as in the proof of Lemma~\ref{lemma-AJ} shows that this lemma is also true for Chow groups with coefficients in $\ZZ $. Moreover, for an element $\alpha \in CH^p_{\ZZ }(S\times X)$ there is the following equivalence: $AJ((\alpha \,_{\CC })_x)=0$ for any closed point $x\in S_{\CC }$ if and only if $AJ(\alpha _x)=0$ for any ($k$-rational) closed point $x\in S$.}
\end{remark}

\begin{remark}
{\rm One can also try to prove an analogue of Lemma~\ref{lemma-AJ} for any, i.e. not necessarily trivial, family of smooth projective varieties $\bcX\to S$, and for an element $\alpha $ in $CH^p(\bcX )$. The implications ``from the right to the left'', in both (i) and (ii), hold true. The implication ``from the left to the right'' in (i) is true as well. The reason for that is that the map $x\mapsto cl(\alpha _x)$, $x\in S_{\CC }$, is a section over $S_{\CC }$ of the local system with fibers $H^{2p}((\bcX _{\CC })_x,\QQ )$. Besides, if
  $$
  AJ((\alpha \,_{\CC })_x)\in (J^{2p-1}((\bcX _{\CC })_x)_{\alg })_{\QQ }
  $$
for any point $x\in S_{\CC }$, then the implication ``from the left to the right'' in (ii) is true by an argument similar to those used in the above proof of (ii).}
\end{remark}

\section{Main result}
\label{mainresult}

The following two propositions can be deduced with the aid of Lemma~\ref{lemma-AJ}. Implicitly they appear in \cite{GreenGriffiths} without a proof (see the first statement in Case 2 on p. 488 in loc.cit.).

\begin{proposition}
\label{corol-1cohom}
Let $X$ and $S$ be two irreducible smooth projective varieties over $k$, $d=\dim(S)$, let $\eta $ be the generic point of $S$, and let $\alpha $ be an element in the group
  $$
  CH^p(S\times X)=CH^{p-d}(S,X)\; .
  $$
Then the following conditions are equivalent:

  \begin{enumerate}

  \item[]{(i)}
  $cl(\alpha _{\eta })=0$;

  \item[]{(ii)}
  the image of the homomorphism $(\alpha\,_{\CC })_*:CH^d(S_{\CC })\to CH^p(X_{\CC })$ is contained in the group $CH^p(X_{\CC})_0$;

  \item[]{(iii)}
  the map $(\alpha _{\CC })_*:H^{2d}(S_{\CC },\QQ )\to H^{2p}(X_{\CC },\QQ )$ is trivial.

 \end{enumerate}

\end{proposition}

\begin{pf}

\medskip

$(i)\Leftrightarrow (ii)$

\medskip

Suppose $(i)$. By  Lemma~\ref{lemma-AJ} (i), $cl((\alpha \,_{\CC })_x)=0$ for any closed point $x$ on $S_{\CC }$. Since, moreover,
$(\alpha \,_{\CC })_x=(\alpha \,_{\CC })_*[x]$, we get
  $$
  cl((\alpha \,_{\CC })_*[x])=0\; ,
  $$
so that $(\alpha \,_{\CC })_*[x]$ is in the group $CH^p(X_{\CC })_0$. As the group $CH^d(S_{\CC })$ is generated by the classes $[x]$, we get $(ii)$.

If we have $(ii)$, then
  $$
  cl((\alpha \,_{\CC })_x)=cl((\alpha \,_{\CC })_*[x])=0
  $$
for any closed point $x\in S_{\CC }$. Applying Lemma~\ref{lemma-AJ} (i), we obtain $(i)$.

\medskip

$(ii)\Leftrightarrow (iii)$

\medskip

Suppose $(ii)$, and take any closed point $x$ in $S_{\CC }$. As the action of correspondences on Chow groups commutes with the action on the cohomology via the cycle class map, and the image of the action of the correspondence $\alpha _{\CC }$ is contained in the group $CH^p(X_{\CC })_0$, we get
  $$
  (\alpha _{\CC })_*(cl([x]))=0\; .
  $$
Since the classes $cl([x])$ generate the $1$-dimensional $\QQ $-vector space $H^{2d}(S_{\CC },{\QQ })$, we get $(iii)$.

Suppose $(iii)$, and take any element $\beta $ in $CH^d(S_{\CC })$. Again, as the action of correspondences on Chow groups commutes with the action on the cohomology, and the map
  $$
  (\alpha _{\CC })_*:H^{2d}(S_{\CC },\QQ )\lra
  H^{2p}(X_{\CC },\QQ )
  $$
is zero, we see that
  $$
  cl((\alpha _{\CC })_*(\beta ))=0\; .
  $$
This gives $(ii)$.
\end{pf}

\begin{proposition}
\label{corol-2cohom}
Let again $X$ and $S$ be two irreducible smooth projective varieties over $k$, $d=\dim(S)$, $\eta $ the generic point of $S$, and let $\alpha $ be an element in the group $CH^p(S\times X)=CH^{p-d}(S,X)$. Assume one of the three equivalent conditions in Proposition~\ref{corol-1cohom}. Then the following conditions are equivalent:

  \begin{enumerate}

  \item[]{(i)}
  $AJ(\alpha _{\eta })=0$;

  \item[]{(ii)}
  the image of the homomorphism $(\alpha\,_{\CC })_*:CH^d(S_{\CC })\to CH^p(X_{\CC })$ is in $T^p(X_{\CC })$;

  \item[]{(iii)}
  the map $(\alpha _{\CC })_*:H^{2d-1}(S_{\CC },\QQ )\to H^{2p-1}(X_{\CC },\QQ )$ is zero, and there exists a closed point $x\in S_{\CC }$, such that $AJ((\alpha\,_{\CC })_x)=0$.

  \end{enumerate}

\end{proposition}

\begin{pf}

\medskip

$(i)\Leftrightarrow (ii)$

\medskip

Suppose $(i)$. By Lemma~\ref{lemma-AJ}$(ii)$, for any closed point $x\in S_{\CC }$ we have that
$AJ((\alpha \,_{\CC })_x)=0$. Since
  $$
  (\alpha \,_{\CC })_*[x]=(\alpha \,_{\CC })_x
  $$
and $AJ((\alpha \,_{\CC })_x)=0$, we obtain
  $$
  AJ((\alpha \,_{\CC })_*[x])=0\; .
  $$
Then, by the definition of the group $T^p(X_{\CC })$, the element $(\alpha\,_{\CC })_*[x]$ sits in $T^p(X_{\CC })$.
Since the group $CH^d(S_{\CC })$ is generated by the classes $[x]$, we arrive to $(ii)$.

Now suppose $(ii)$. Then
  $$
  AJ((\alpha _{\CC })_*[x])=0
  $$
for any closed point $x\in S_{\CC }$. Since $(\alpha _{\CC })_*[x]=(\alpha _{\CC })_x$, we see that $AJ((\alpha _{\CC })_x)=0$ for any closed point $x\in S_{\CC }$. Applying Lemma~\ref{lemma-AJ} (ii), we get $(i)$.

\medskip

$(ii)\Leftrightarrow (iii)$

\medskip

Assuming $(ii)$, for any closed point $x\in S_{\CC }$, one has
  $$
  AJ((\alpha \,_{\CC })_x)=AJ((\alpha \,_{\CC})_*[x])=0\; .
  $$
So, to deduce (iii) we need only to show that the map
  $$
  (\alpha _{\CC })_*:H^{2d-1}(S_{\CC },\QQ )\lra
  H^{2p-1}(X_{\CC },\QQ )
  $$
is zero. Let
  $$
  \gamma \in CH^d(S_{\CC })_0
  $$
be a zero-cycle class of degree zero. We have that
  $$
  AJ((\alpha \,_{\CC })_*(\gamma ))=(\alpha \,_{\CC })_*(AJ(\gamma ))\; .
  $$
As the element $(\alpha\,_{\CC })_*(\gamma )$ sits in $T^p(X_{\CC })$, the element $AJ((\alpha\,_{\CC })_*(\gamma ))$ vanishes, so that
  $$
  (\alpha \,_{\CC })_*(AJ(\gamma ))=0\; .
  $$

Since the Albanese map
  $$
  AJ : CH^d(S_{\CC })_0\to J^{2d-1}(S_{\CC })_\QQ =
  \Alb(S_{\CC })_{\QQ }
  $$
is surjective and $(\alpha \,_{\CC })_*(AJ(\gamma ))=0$, the map
  $$
  (\alpha \,_{\CC })_* : J^{2d-1}(S_{\CC })_{\QQ }\to J^{2p-1}(X_{\CC })_{\QQ }
  $$
is zero.

Let $\beta \in CH^p_{\ZZ }(S\times X)$ be an integral cycle class, such that $\beta =n\cdot \alpha $ for some non-zero integer $n$. The map
  $$
  (\beta \,_{\CC })_*:J^{2d-1}(S_{\CC })\to J^{2p-1}(X_{\CC })
  $$
satisfies
  $$
  (\beta \,_{\CC })_*\otimes_{\ZZ }\QQ =
  n\cdot (\alpha \,_{\CC })_*\; .
  $$
Therefore,
  $$
  (\beta _{\CC })_*\otimes_{\ZZ }\QQ =0\; .
  $$
It means that the image of $(\beta \,_{\CC })_*$ is contained in the torsion of the group $J^{2p-1}(X_{\CC })$. Since $(\beta _{\CC })_*$ is a continuous map between complex compact tori, the image of $(\beta _{\CC })_*$ is compact and connected. As the subgroup $(\beta _{\CC })_*(J^{2d-1}(S_{\CC }))$ is compact, connected and torsion in the compact complex torus $J^{2p-1}(X_{\CC })$, it is trivial. Consequently, the map
  $$
  (\beta _{\CC })_*:J^{2d-1}(S_{\CC })\to
  J^{2p-1}(X_{\CC })
  $$
vanishes, and the corresponding map on first homology groups
  $$
  H_1(J^{2d-1}(S_{\CC }),\QQ )\to
  H_1(J^{2p-1}(X_{\CC }),\QQ )\; ,
  $$
induced by $(\beta _{\CC })_*$, vanishes too.

By the construction of the intermediate Jacobians,
  $$
  H_1(J^{2d-1}(S_{\CC }),\QQ ) = H^{2d-1}(S_{\CC },\QQ )
  $$
and
  $$
  H_1(J^{2p-1}(X_{\CC }),\QQ ) = H^{2p-1}(X_{\CC },\QQ )\; ,
  $$
and the above action on first homology groups agrees with the map
  $$
  (\beta _{\CC })_*:H^{2d-1}(S_{\CC },\QQ )\to
  H^{2p-1}(X_{\CC },\QQ )\; ,
  $$
so that it vanishes. Since $\alpha =\frac{1}{n}\cdot \beta $
in $CH^p(S\times X)$, the map
  $$
  (\alpha \, _{\CC })_*:H^{2d-1}(S_{\CC },\QQ )\to
  H^{2p-1}(X_{\CC },\QQ )
  $$
vanishes. This gives $(iii)$.

Now suppose $(iii)$. By the construction of the groups $J^{2d-1}(S_{\CC })_{\QQ }$ and $J^{2p-1}(X_{\CC })_{\QQ }$, the vanishing of the map $(\alpha _{\CC })_*:H^{2d-1}(S_{\CC },\QQ )\to H^{2p-1}(X_{\CC },\QQ )$ implies vanishing of the map
  $$
  (\alpha _{\CC })_*:J^{2d-1}(S_{\CC })_{\QQ }\lra
  J^{2p-1}(X_{\CC })_{\QQ }\; .
  $$
Therefore, for any zero-cycle $\beta \in CH^d(S_{\CC })_0$, we have that
  $$
  (\alpha _{\CC })_*(AJ(\beta ))=0\; .
  $$
Since $AJ((\alpha _{\CC })_*(\beta))=(\alpha _{\CC })_*(AJ(\beta))$, the element $AJ((\alpha _{\CC })_*(\beta))$ vanishes too. Therefore, the image of the map
  $$
  (\alpha _{\CC })_*:CH^d(S_{\CC })_0\lra CH^p(X_{\CC })_0
  $$
is contained in the group $T^p(X_{\CC })$.

Now take the closed point $x\in S_{\CC }$ mentioned in $(iii)$. As $AJ((\alpha _{\CC })_x)=0$ and $(\alpha _{\CC })_*[x]=(\alpha\,_{\CC })_x$, we have that
  $$
  (\alpha _{\CC })_*[x]\in T^p(X_{\CC })\; .
  $$
Since the whole group $CH^d(S_{\CC })$ is generated by the subgroup $CH^d(S_{\CC })_0$ and the class $[x]$, we obtain $(ii)$.
\end{pf}

\bigskip

Now we want a sufficient condition for $\alpha _{\eta }$ to be a non-zero class. It is equivalent to non-vanishing of some Hodge-theoretic invariants constructed in \cite{GreenGriffiths}. More precisely, the condition in Theorem~\ref{non-triviality} below is equivalent to the vanishing of invariants denoted by $[\mathscr Z]_m$ and defined on p. 483 in loc.cit.

\begin{theorem}
\label{non-triviality}
Let $X$ and $S$ be smooth projective irreducible varieties over $k$, let $d$ be the dimension of $S$, $\eta $ the generic point of $S$, and let $\alpha $ be an element in the group
  $$
  CH^p(S\times X)=CH^{p-d}(S,X)\; .
  $$
Suppose that there exists $i$, such that the map
  $$
  (\alpha _{\CC })_*:H^{i,d}(S_{\CC })\lra H^{i+p-d,p}(X_{\CC })
  $$
is non-zero. Then the cycle class $\alpha _{\eta }$ is non-zero in $CH^p(X_{k(S)})$.
\end{theorem}

\begin{pf}
Suppose $\alpha _{\eta }=0$ in $CH^p(X_{k(S)})$. It is equivalent to say that $\alpha $ is balanced on the left, see \cite{KMP}. By Lemma~\ref{lemma-suppleft}, there exists an equidimensional smooth projective variety $Z$ over $k$ with $\dim (Z)<d$, such that the morphism of motives
  $$
  \alpha _*:M(S)\lra M(X)\otimes \Le ^{\otimes (d-p)}
  $$
is equal to a composition
  $$
  M(S)\to M(Z)\to M(X)\otimes \Le ^{\otimes (d-p)}\; .
  $$
Therefore, the non-zero map
  $$
  (\alpha _{\CC })_*:H^{i,e}(S_{\CC })\lra H^{i+p-e,p}(X_{\CC })
  $$
factors through the group $H^{i,e}(Z_{\CC })$. As $\dim (Z)<d$, the latter group is trivial, which is in contradiction with the non-triviality of the map $\alpha _*$.
\end{pf}

\bigskip

Now we are almost ready to give a method of constructing non-trivial elements in the Abel-Jacobi kernels. By \cite{Murre}, for an irreducible smooth projective variety $S$ of dimension $d$ over $k$, there exist the Picard and the Albanese projectors of $S$ generating, respectively, the Picard motive $M^1(S)$ and its dual, the Albanese motive $M^{2d-1}(S)$, of the variety $S$. Denote by
  $$
  \pi_{2d-1}\in CH^0(S,S)
  $$
the projector that corresponds to the motive $M^{2d-1}(S)$. Let also
  $$
  \pi_{2d}\in CH^0(S,S)
  $$
be a projector determined by a fixed closed point on $S$ and corresponding to a direct summand $\Le^{\otimes d}$ in $M(S)$.

Put
  $$
  \tau = [\Delta ]-\pi _{2d-1}-\pi _{2d}\; ,
  $$
where $[\Delta ]$ is the class of the diagonal in the group $CH^0(S,S)$. Then $\tau _{\CC }$ acts as zero on the groups $H^{2d}(S_{\CC },\QQ )$ and $H^{2d-1}(S_{\CC },\QQ )$, and identically on the groups $H^i(S_{\CC },\QQ )$ with $i<2d-1$.

If now $X$ is another irreducible smooth projective variety over $k$ and $\alpha $ is a correspondence of degree $d-p$ from $S$ to $X$, i.e. an element in the group $CH^p(S\times X)=CH^{p-d}(S,X)$, let
  $$
  \alpha '=\alpha \circ \tau \in CH^{p-d}(S,X)
  $$
be a correspondence modified by the above correspondence $\tau $, and if $x$ is a closed point on $S$, then let
  $$
  \alpha '(x) =\alpha '-[S]\times \alpha '_x\in CH^{p-d}(S,X)
  $$

\medskip

be a correspondence of the same degree $p-d$ from $S$ to $X$ normalized by the point $x$.

Our main result in this paper is as follows:

\medskip

\begin{theorem}
\label{maintheorem}
Let $X$ and $S$ be two irreducible smooth projective varieties over $k$, $d=\dim(S)$, and let $\alpha $ be an element in
  $$
  CH^p(S\times X)=CH^{p-d}(S,X)\; .
  $$
Suppose there exists $i\leq d-2$, such that the induced map
  $$
  (\alpha _{\CC })_*:H^{i,d}(S_{\CC })\lra H^{i+p-d,p}(X_{\CC })
  $$
is non-zero. Then, for any closed point $x$ in $S$, the specialization
  $$
  \alpha '(x)_{\eta }\in CH^p(X_{k(S)})
  $$
of the cycle class $\alpha '(x)$ at the generic point $\eta $ is a non-zero element in the Abel-Jacobi kernel $T^p(X_{k(S)})$.
\end{theorem}

\begin{pf}
The correspondence $([S]\times \alpha '_x)_{\CC }=[S_{\CC }]\times (\alpha '_{\CC })_x$ acts as zero on all the groups $H^i(S_{\CC },\QQ )$ with $i\ne 2d$ by a dimension reason. By the construction of $\tau $, the modified correspondence $\alpha '_{\CC }$ acts as zero on the group $H^{2d}(S_{\CC },\QQ )$. As the action of correspondences on Chow groups is compatible with their action on cohomology,
  $$
  (\alpha '_{\CC })_*(CH^d(S_{\CC }))\subset
  CH^p(X_{\CC })_0\; .
  $$
In particular,
  $$
  cl((\alpha '_{\CC })_x)=cl((\alpha '_{\CC })_*[x])=0\; .
  $$
Therefore, for any element $\gamma \in H^{2d}(S_{\CC },\QQ )$ we have that
  $$
  (([S]\times \alpha '_x)_{\CC })_*(\gamma )=
  cl((\alpha '_{\CC })_x)\cdot \deg (\gamma )=
  0\cdot \deg (\gamma )=0\; ,
  $$
where
  $$
  \deg : H^{2d}(S_{\CC },\QQ )\lra \QQ
  $$
is the orientation isomorphism. Thus, the correspondence $([S]\times \alpha '_x)_{\CC }$ acts as zero on all the cohomology groups $H^i(S_{\CC },\QQ )$.

The correspondence $\alpha '_{\CC }$ acts as zero not only on the group $H^{2d}(S_{\CC },\QQ )$, but also on the group group $H^{2d-1}(S_{\CC },\QQ)$. Consequently, the correspondence $\alpha '(x)_{\CC }$ acts as zero on the groups $H^{2d}(S_{\CC },\QQ)$ and $H^{2d-1}(S_{\CC },\QQ )$.

Next, since
  $$
  \alpha '(x)_*([x])=
  \alpha '_*([x])-([S]\times \alpha '_x)_*([x])=
  \alpha '_x-\alpha '_x=0\; ,
  $$
we have that
  $$
  (\alpha '(x)_{\CC })_x=(\alpha '(x)_*([x]))_{\CC }=0\; .
  $$

Then we see that the correspondence $(\alpha '(x))$ satisfies the conditions of Proposition~\ref{corol-1cohom} (iii) and Proposition~\ref{corol-2cohom} (iii), so that $(\alpha '(x))_{\eta }$ is in $T^p(X_{k(S)})$.

Moreover, the action of the correspondence $\alpha '(x)_{\CC }$  coincides with the action of the correspondence $\alpha _{\CC }$ on the groups $H^i(S_{\CC },\QQ )$ with $i\leq 2d-2$. Since, by the assumption of the theorem, $\alpha _{\CC }$ acts non-trivially on $H^i(S_{\CC },\QQ )$ for some $i\leq 2d-2$, we see that $\alpha '(x)_{\CC }$ acts non-trivially on $H^i(S_{\CC },\QQ )$ too. By Theorem~\ref{non-triviality}, the element $\alpha '(x)_{\eta }$ is non-trivial in $T^p(X_{k(S)})$.
\end{pf}

\begin{remark}
\label{lower}
{\rm Certainly, we can always reduce the dimension of $S$ from
Theorem~\ref{maintheorem} with the aid of the Lefschetz hyperplane section theorem. Indeed, let $X$, $S$ and $\alpha $ be as in Theorem \ref{maintheorem}. Let $T$ be a general multiple hyperplane section of $S$ of codimension $i$, whose generic point will be denoted by $\xi $. Let $f:T\hra S$ be the corresponding closed embedding, and let
  $$
  \beta =(f\times \id _X)^*(\alpha )\in CH^{p-d+i}(T,X)
  $$
be the pull-back of the correspondence $\alpha $ to $T\times X$. The composition
  $$
  H^{0,d-i}(T_{\CC })\stackrel{(f_{\CC })_*}{\lra }H^{i,d}(S_{\CC })\stackrel{(\alpha _{\CC })_*}\lra H^{i+p-d,p}(X_{\CC })
  $$
is equal to the action
  $$
  (\beta _{\CC })_* : H^{0,d-i}(T_{\CC })\lra
  H^{i+p-d,p}(X_{\CC })\; .
  $$
As $\dim(T)=d-i$, by the Lefschetz hyperplane section theorem, the homomorphism $(f_{\CC })_*$ is surjective. Therefore, the non-vanishing of $(\alpha _{\CC })_*$ implies non-vanishing of $(\beta _{\CC })_*$. We see that $\beta $ satisfies the condition of Theorem \ref{maintheorem} too. As in the theorem, for any closed point $x$ on $T$, we look at the cycle class
  $$
  \beta '(x) = \beta ' - [T]\times \beta '_x
  $$
in the group $CH^p(T\times X)$. Then the specialization at the generic point,
  $$
  \beta '(x)_{\xi }\in CH^p(X_{k(T)})
  $$
is a non-trivial element in the Abel-Jacobi kernel $T^p(X_{k(T)})$.
}
\end{remark}

\begin{remark}
\label{remark-component}
{\rm Under the notation of Theorem~\ref{maintheorem}, consider the class
  $$
  cl(\alpha _{\CC })=\sum _{r+s=2p} cl(\alpha _{\CC })_{r,s}\in
  H^{2p}((S\times X)_{\CC },{\CC })\cong
  \bigoplus_{r+s=2p}
  H^s(S_{\CC },\CC )\otimes _{\CC }H^r(X_{\CC },\CC )\; .
  $$
The Hodge theory implies that the correspondence $\alpha $ satisfies the condition of Theorem~\ref{maintheorem} if and only if the image of the class $cl(\alpha _{\CC })_{2p+i-d,d-i}$ under the projection
  $$
  H^{d-i}(S_{\CC },\CC )\otimes _{\CC }
  H^{2p+i-d}(X_{\CC },\CC )\to
  H^{d-i,0}(S_{\CC })\otimes_{\CC }H^{i+p-d,p}(X_{\CC })
  $$
does not vanish.
}
\end{remark}

\section{Some examples}
\label{someexamples}

In this section we outline how to apply Theorem \ref{maintheorem} to construct algebraic cycle classes in the Abel-Jacobi kernels for $K3$-surfaces, which is classic, and then for threefolds and higher-dimensional varieties with known the Lefschetz conjecture.

\medskip

Let first $X$ be a smooth projective surface over $k$, let $S=X$ and let $p=2$. Then a correspondence $\alpha \in CH^2(X\times X)=CH^0(X,X)$ satisfies the condition of Theorem~\ref{maintheorem} if and only if the map
  $$
  \alpha _*:H^{0,2}(X)\lra H^{0,2}(X)
  $$
is non-zero. Suppose, for example, $X$ is a $K3$ surface with an algebraic symplectomorphism $f:X\to X$ on it. Without loss of generality one can assume that $k$ is big enough so that $s$ is defined over $k$. Then take $\alpha $ to be the class of the graph $\Gamma _f$ of the symplectomorphism $f$.

\medskip

Now let $X$ be a smooth projective irreducible variety over $k$ of dimension $d$, such that, for some $p\geq 2$, the group $H^{p,0}(X_{\CC })$ is non-trivial, and the inverse of the Lefschetz operator
  $$
  L : H^p(X_{\CC },\QQ )\lra H^{2d-p}(X_{\CC },\QQ )
  $$
is represented by an algebraic correspondence
  $$
  \alpha \in CH^p(X\times X)=CH^{p-d}(X,X)\; .
  $$

Let $Y$ be an intersection of general $d-p$ hyperplane sections of $X$, and let
  $$
  f : Y\hra X
  $$
be the corresponding closed embedding. Then one has a morphism of Chow motives over the field $k$,
  $$
  \beta =\alpha \circ M(f)^{\rm t}:M(Y)\lra M(X)\; ,
  $$
where $M(f)^{\rm t}$ is the transposition of the graph of the embedding $f$. Notice that $\beta $ is an element in the group $CH^p(Y\times X)=CH^0(Y,X)$.

Under the above notation, the elements $\alpha $ and $\beta $ satisfy the condition of Theorem~\ref{maintheorem} because, by the Lefschetz theorem, the maps
  $$
  (\alpha _{\CC })_*:H^{d-p,d}(X_{\CC })\lra
  H^{0,p}(X_{\CC })\; ,
  $$
  $$
  (\beta _{\CC })_*:H^{0,p}(S_{\CC })\lra H^{0,p}(X_{\CC })
  $$
are both isomorphisms.

Let $\eta $ be the generic point of $X$, and let $\xi $ be the generic point of $Y$. Take any two closed points $x\in X$ and $y\in Y$. Then we obtain two non-zero elements
  $$
  \alpha '(x)_{\eta }=\alpha '-[X]\times \alpha '_x\in T^p(X_{k(X)})
  $$
and
  $$
  \beta '(y)_{\xi }=\beta ' - [Y]\times \beta '_y\in T^p(X_{k(Y)})\; .
  $$

\bigskip

Let us consider in some more detail the case when $d=3$ and $p=2$. The existence of $\alpha $ as above is equivalent to say that the Lefschetz conjecture holds for $X$, which implies the existence of a K\"unneth decomposition for the homological motive of $X$. Assuming finite-dimensionality of Chow motive $M(X)$ or, equivalently, that the homological realization functor is conservative, we obtain a Chow-K\"unneth decomposition for the motive of $X$, see \cite{GP2}:
  $$
  M(X)\cong\oplus_{i=0}^6 M^i(X)\; .
  $$
As $Y$ is a surface, it also has a Chow-K\"unneth decomposition for its motive $M(Y)$. The morphism
  $$
  \beta : M(Y)\lra M(X)
  $$
induces the morphisms
  $$
  \beta : M^i(Y)\to M^i(X),\quad i\geq 0\;,
  $$
and $M^2(X)$ is a direct summand of the middle motive $M^2(Y)$. According to the conjectural formulas for the Bloch-Beilinson filtration, for any field extension $k\subset K$, we have that
  $$
  T^2(X_K)=CH^2(M^2(X)_K)
  $$
(possibly, this also follows from finite-dimensionality).

Thus, the big conjectures applied to $X$ imply that the map
  $$
  (\beta _K)_*:T^2(Y_K)\lra T^2(X_K)
  $$
is surjective.

One may wish to ask the following question: is it possible to prove the above statement unconditionally?

Notice that, by construction, the element $\beta '(y)_{\xi }$ belongs to the image of the map
  $$
  (\beta _{k(Y)})_*:T^2(Y_{k(Y)})\lra T^2(X_{k(Y)})\; .
  $$

A more particular question: is it possible to prove unconditionally that the element $\alpha '(x)_{\eta }$ belongs to the image of the map
  $$
  (\beta _{k(X)})_*:T^2(Y_{k(X)})\lra T^2(X_{k(X)})\, ?
  $$

The above questions can be tested on a wide range of examples of threefolds that satisfy the Lefschetz conjecture, \cite{Tankeev}, and whose motive is finite-dimensional, \cite{GG}.

\section{Self-intersection of the Poincar\'e bundle}
\label{mainexample}

In this section we will demonstrate another method how to produce non-trivial elements in the Abel-Jacobi kernels for high-dimensional varieties, with the aid of Theorem \ref{non-triviality}.

Let $X$ be an irreducible smooth projective variety over $k$, and let $A$ be the Picard variety $\Pic ^0(X)$ of the variety $X$. Fix a closed point $x_0\in X$ and consider the associated Poincar\'e line bundle $\bcP $ over $A\times X$. Recall that $\bcP $ is the unique line bundle on $A\times X$, such that for any point $a\in A$ the isomorphism class of the line bundle on $X$
  $$
  \bcP |_{a\times \, X}
  $$
is equal to $a$, and
  $$
  \bcP |_{\, A\times x_0}\cong{\mathscr O}_A\; .
  $$
Let
  $$
  \beta \in CH^1(A\times X)=CH^{1-g}(A,X)
  $$
be the first Chern class of the bundle $\bcP $, where
  $$
  g=\dim (A)
  $$
is the dimension of $A$. For each index $i$ the correspondence $\beta $ induces a map on cohomology:
  $$
  (\beta _{\CC })_{*,i}:H^i(A_{\CC },\QQ )\lra
  H^{i+2-2g}(X_{\CC },\QQ )\; .
  $$

\begin{lemma}
\label{lemma-isom}
The homomorphism $(\beta _{\CC })_{*,i}$ is an isomorphism for $i=2g-1$, and vanishes otherwise.
\end{lemma}

\begin{pf}
The map $(\beta _{\CC })_{*,2g-1}$ is an isomorphism because the correspondence $\beta $ induces an isomorphism between the Albanese variety $\Alb (A)$ and the Picard variety $\Pic ^0(X)$. Thus, we have only to show that $(\beta _{\CC })_{*,i}$ vanishes when $i\neq 2g-1$.

As the source and the target of a non-trivial homomorphism need to be non-zero, we see that $(\beta _{\CC })_{*,i}$ can be non-trivial only when
$2g-2\leq i\leq 2g$. We have already shown that $cl((\beta _{\CC })_a)=0$ for any point $a$ in $A_{\CC }$. Hence, the map $(\beta _{\CC })_{*,2g}$ vanishes by Proposition \ref{corol-1cohom}, and we have only to deduce that $(\beta _{\CC })_{*,2g-2}$ is zero.

If $i=2g-2$, the group in the target, i.e. $H^0(X_{\CC },\QQ )$, is canonically isomorphic to $\QQ $, and for any closed point
  $$
  i_x : x\lra X_{\CC }
  $$
the pull-back
  $$
  i_x^* : H^0(X_{\CC },\QQ )\lra H^0(x,\QQ )
  $$
is an isomorphism, so that $(\beta _{\CC })_{*,2g-2}$ vanishes if and only if the composition $i_x^*\circ (\beta _{\CC })_{*,2g-2}$ vanishes.

Let $\gamma $ be an arbitrary element in the group $H^{2g-2}(A_{\CC },\QQ )$. By the definition of the action of correspondences, we have that
  $$
  i_x^*(\beta _{\CC })_*(\gamma )=i_x^*(p_X)_*(p_A^*\gamma \cdot cl(\beta _{\CC }))\; .
  $$
Now look at the Cartesian square
  $$
  \diagram
  A\times x\ar[dd]_-{p_x} \ar[rr]^-{\id _A\times i_x} & &
  A\times X \ar[dd]^-{p_X} \\ \\
  x \ar[rr]^-{i_x} & & X
  \enddiagram
  $$
where $p_X$ and $p_x$ are the natural projections. Since both $i_x$ and $\id _A\times i_x$ are regular embeddings of the same codimension, we have that
  $$
  i_x^*(p_X)_*(p_A^*\gamma \cdot cl(\beta \,_{\CC }))=
  (p_x)_*(\id _A\times i_x)^*(p_A^*\gamma \cdot cl(\beta \,_{\CC }))\; ,
  $$
see \cite{Fulton}, Theorem 6.2. Pull-backs commute with products in cohomology, so that
  $$
  (p_x)_*(\id _A\times i_x)^*(p_A^*\gamma \cdot cl(\beta _{\CC }))=
  (p_x)_*((\id_A\times i_x)^*(p_A^*\gamma)\cdot(\id _A\times i_x)^*(cl(\beta _{\CC })))\; .
  $$
The cycle class $(\id _A\times i_x)^*(p_A^*\gamma)$ coincides with $\gamma$, under the identification of $A\times x$ with $A$. Therefore,
  $$
  (p_x)_*((\id _A\times i_x)^*(p_A^*\gamma )\cdot (\id _A\times i_x)^*(cl(\beta \,_{\CC })))=
  (p_x)_*(\gamma \cdot (\id _A\times i_x)^*(cl(\beta \,_{\CC })))\; .
  $$
Since the cycle class map commutes with pull-backs, we have:
  $$
  (p_x)_*(\gamma \cdot (\id _A\times i_x)^*(cl(\beta\,_{\CC })))=
  (p_x)_*(\gamma \cdot cl((\id _A\times i_x)^*\beta \,_{\CC }))\; .
  $$
The element $(\id _A\times i_x)^*(\beta _{\CC })$ is equal to $(\beta ^{\rm t}_{\CC })_x$, where $\beta ^{\rm t}$ is the transpose of the correspondence $\beta $. In particular, $\beta ^{\rm t}$ is an element in the Chow group $CH^{1-\dim (X)}(X,A)$.

As a result, we obtain:
  $$
  i_x^*(\beta _{\CC })_*(\gamma )=
  (p_x)_*(\gamma \cdot cl((\beta _{\CC }^{\rm t})_x))\; .
  $$

Next, $(\beta _{\CC }^{\rm t})_{x_0}$ is the first Chern class of the line bundle $\bcP |_{A\times x_0}$, where $x_0\in X$ is the closed point that has been fixed above for the choice of the Poincar\'e line bundle. Then, by the definition of $\bcP $, the line bundle $\bcP |_{A\times x_0}$ is trivial, whence $(\beta _{\CC }^{\rm t})_{x_0}=0$. In particular,
  $$
  cl((\beta \,_{\CC }^{\rm t})_{x_0})=0\; .
  $$
As the cycle classes $(\beta _{\CC }^{\rm t})_x$ and $(\beta _{\CC }^{\rm t})_{x_0}$ are algebraically equivalent,
  $$
  cl((\beta _{\CC }^{\rm t})_{x})=0
  $$
too. Since $i_x^*(\beta _{\CC })_*(\gamma )=(p_x)_*(\gamma \cdot cl((\beta _{\CC }^{\rm t})_{x}))$, one has:
  $$
   i_x^*(\beta _{\CC })_*(\gamma )=0\; ,
   $$
for any $\gamma $ in $H^{2g-2}(A_{\CC },\QQ )$. This implies the vanishing
of $(\beta _{\CC })_{*,2g-2}$.
\end{pf}

\bigskip

Let now
  $$
  \alpha =\beta ^p\in CH^p(A\times X)
  $$
be the $p$-fold self-intersection of the class $\beta $,
where $p\geq 2$. Let $S=A$ and let now
  $$
  \eta =\Spec (k(A))
  $$
be the generic point of the abelian variety $A$. Then one has a restriction
  $$
  \alpha _{\eta }\in CH^p(X_{k(A)})
  $$
of the correspondence $\alpha $ on the generic point from the left.

\medskip

\begin{theorem}
\label{oneforms}
Under the above notation, the element $\alpha _{\eta }$ belongs to the Abel-Jacobi kernel $T^p(X_{k(A)})$. Moreover, if there exist $p$ differential forms
  $$
  \omega _1,\ldots ,\omega _p\in
  H^0(X_{\CC },\Omega ^1_{X_{\CC }})=H^{1,0}(X_{\CC })\; ,
  $$
such that
  $$
  \omega _1\wedge \ldots \wedge \omega _p\neq 0\; ,
  $$
then
  $$
  \alpha _{\eta }\neq 0
  $$
in $T^p(X_{k(A)})$.
\end{theorem}

\begin{pf}
By the construction of the cycle class $\beta $, for any closed point $a\in A_{\CC }$ the element $(\beta \, _{\CC })_a$ is equal, in the group $CH^1(X_{\CC })$, to the first Chern class of the line bundle  $\bcP|_{X\times a}$. Hence, the class $cl((\beta \,_{\CC })_a)$ is equal to the first Chern class of $\bcP |_{X\times a}$ with value in the second Betti cohomology group $H^2(X_{\CC },\QQ )$. By the definition of the Poincar\'e line bundle, for any $a\in \Pic^0(X_{\CC })$ the isomorphism class of $\bcP |_{X\times a}$ is exactly $a$. Remind also that elements in $\Pic^0(X_{\CC })$ are isomorphism classes of line bundles with a trivial first Chern class with value in the second Betti cohomology group $H^2(X_{\CC },\QQ )$. Therefore, we have that
  $$
  cl((\beta \,_{\CC })_a)=0\; ,
  $$
for any closed point $a$ in $A_{\CC }$.

Let
  $$
  \beta _{\eta }\in CH^1(X_{k(A)})
  $$
be the restriction if $\beta $ at the generic point of $A$. By Lemma~\ref{lemma-AJ}$(i)$, we have that
  $$
  cl(\beta _{\eta })=0\; .
  $$

Fix an embedding of the function field $k(A)$ in to $\CC $ over $k$. Then it makes sense to speak about scalar extensions like $(\alpha _{\eta })_{\CC }$ or $(\beta _{\eta })_{\CC }$. Since
  $$
  cl((\alpha _{\eta })_{\CC })=cl((\beta _{\eta })_{\CC }^p)\; ,
  $$
where $(\beta _{\eta })_{\CC }^p$ is the $p$-fold self-intersection of the cycle class $(\beta _{\eta })_{\CC }$ in $X_{\CC }$, we have that
  $$
  cl(\alpha _{\eta })=0\; .
  $$
Appling $AJ$ to the cycle class $(\alpha _{\eta })_{\CC }$ we obtain:
  $$
  AJ((\alpha _{\eta })_{\CC })=
  AJ((\beta _{\eta })_{\CC }^{p-1}\cdot
  (\beta _{\eta })_{\CC })=
  cl((\beta _{\eta })_{\CC }^{p-1})\cdot
  AJ((\beta _{\eta })_{\CC })
  $$
by Proposition 9.23 in \cite{Voisin}. Since $p\geq 2$, and $cl(\beta _{\eta })=0$, we get
  $$
  cl((\beta _{\eta })_{\CC }^{p-1})=0\; .
  $$
It follows that
  $$
  AJ((\alpha _{\eta })_{\CC })=0\; .
  $$
Thus,
  $$
  \alpha _{\eta }\in T^p(X)\; .
  $$

The rest of the proof of Theorem \ref{oneforms} is devoted to showing that $\alpha $ satisfies the condition of Theorem~\ref{non-triviality}, which will allow us to prove that $\alpha _{\eta }$ is non-zero.

\medskip

With this aim we will find the shape of the K\"unneth components of the class $cl(\beta _{\CC })$ in the cohomology groups of the product $A\times X$. Taking products in cohomology gives rise to the K\"unneth components for the cycle class $\alpha _{\CC }=\beta _{\CC }^p$. This will allow to describes explicitly the action of $\alpha _{\CC }$ on $H^{2g-p}(A_{\CC },\QQ )$. Then, using the fact that $A$ is an abelian variety, we will construct an explicit element in $H^{g,g-p}(A_{\CC })$ with a non-trivial image under the map $(\alpha _{\CC })_*$. Finally, we will use Theorem~\ref{non-triviality} in order to show that $\alpha _{\eta }$ is non-zero.

\medskip

So, let us implement this program. Consider the K\"unneth components of $cl(\beta \,_{\CC })$ in the cohomology group $H^*(A_{\CC },\QQ)\otimes_\QQ H^*(X_{\CC },\QQ)$:
  $$
  cl(\beta \,_{\CC })=\sum_{i+j=2}cl(\beta \,_{\CC })_{i,j}\in \oplus_{i+j=2}
  H^i(A_{\CC },\QQ)\otimes_{\QQ }H^j(X_{\CC },\QQ)\; .
  $$
The Poincar\'e duality establishes isomorphisms for all $i$:
  $$
  H^i(A_{\CC },\QQ )\cong H^{2g-i}(A_{\CC },\QQ)^\vee \; ,
  $$
where by $V^\vee$ we denote the dual vector space to a finite-dimensional vector space $V$. Therefore, the tensor product
  $$
  H^i(A_{\CC },\QQ )\otimes_{\QQ }H^j(X_{\CC },\QQ)
  $$
is canonically isomorphic to the space of $\QQ $-linear operators
  $$
  \Hom _{\QQ }(H^{2g-i}(A_{\CC },\QQ ),H^j(X_{\CC },\QQ ))\; ,
  $$
so that one can consider elements $cl(\beta _{\CC })_{i,j}$ in tensor products as operators. With this identification, for any element
  $$
  \gamma \in H^{2g-i}(A_{\CC },\QQ )
  $$
we have that
  $$
  (\beta _{\CC })_*\gamma =cl(\beta _{\CC })_{i,j}(\gamma )\; .
  $$
In the above notation, we have that
  $$
  (\beta _{\CC })_{*,i}=cl(\beta _{\CC })_{2g-i,2-2g+i}\; .
  $$
Using the known vanishing of $(\beta _{\CC })_{*,i}$, see Lemma \ref{lemma-isom}, we deduce that
  $$
  cl(\beta _{\CC })=
  cl(\beta _{\CC })_{1,1}\in H^1(A_{\CC },\QQ )\otimes_{\QQ }
  H^1(X_{\CC },\QQ )\; .
  $$
By the construction of $\alpha $,
  $$
  cl(\alpha \,_{\CC })=cl(\beta \,_{\CC })^p\; .
  $$
And, as we have seen above,
  $$
  cl(\beta _{\CC })\in
  H^1(A_{\CC },\QQ )\otimes_{\QQ }H^1(X_{\CC },\QQ )\; .
  $$
Therefore, it makes sense to compute the K\"unneth components of $p$-th powers of a given element in $H^1(A_{\CC },\QQ)\otimes_\QQ H^1(X_{\CC },\QQ)$.

\medskip

Let
  $$
  \nu : H^1(A_{\CC },\QQ )\otimes_{\QQ }H^1(X_{\CC },\QQ )\to H^{2p}(A_{\CC }\times X_{\CC },\QQ )
  $$
be a map sending
  $$
  \sum _i a_i\otimes x_i\mapsto \left(\sum_i a_i\otimes x_i\right)^p\; ,
  $$
where $a_i\in H^1(A_{\CC },\QQ )$ and $x_i\in H^1(X_{\CC },\QQ )$.

Notice that the product in cohomology groups of $A\times X$ factors into the product in the cohomology groups of $A$ and the product in the cohomology groups of $X$. In other words, the algebra $H^*(A_{\CC }\times X_{\CC },\QQ)$ is isomorphic to the tensor product of the algebras $H^*(A_{\CC },\QQ )$ and
$H^*(X_{\CC },\QQ )$. Another thing is that the product in the first cohomology group is skew-symmetric.

Denote by $\diamond $ the multiplication in the tensor product
  $$
  \wedge ^*H^1(A_{\CC },\QQ )\otimes_{\QQ }
  \wedge ^*H^1(X_{\CC },\QQ )
  $$
of the Grassman algebras. Let
  $$
  \lambda :H^1(A_{\CC },\QQ )\otimes _{\QQ }
  H^1(X_{\CC },\QQ )\to
  \wedge ^pH^1(A_{\CC },\QQ )\otimes_{\QQ }\wedge^p
  H^1(X_{\CC },\QQ )\; ,
  $$
be a map sending
  $$
  \sum _i a_i\otimes x_i\mapsto
  \left(\sum _i a_i\otimes x_i\right)^{\diamond p}\; ,
  $$
where $\diamond $ is the just defined multiplication, and let
  $$
  \mu :\wedge ^p H^1(A_{\CC },\QQ )\otimes_{\QQ }
  \wedge ^p H^1(X_{\CC },\QQ )\to
  H^p(A_{\CC },\QQ )\otimes_{\QQ }H^p(X_{\CC },\QQ )\; ,
  $$
be a map given by the multiplication in cohomology and the natural embedding
  $$
  H^p(A_{\CC },\QQ )\otimes_{\QQ }H^p(X_{\CC },\QQ )\to
  H^{2p}(A_{\CC }\times X_{\CC },\QQ )\; .
  $$
Then
  $$
  \nu =\mu \circ \lambda \; .
  $$
We need to describe the maps $\lambda $ and $\mu $ in terms of operators between cohomology groups.

Since the wedge powers of dual vector spaces are still dual, there is a canonical isomorphism
  $$
  \wedge ^pH^1(A_{\CC },\QQ )\cong
  \wedge ^pH^{2g-1}(A_{\CC },\QQ )^{\vee }\; .
  $$
Therefore, the target of $\lambda $ and the source of $\mu $ is the space of operators
  $$
  \Hom _{\QQ }(\wedge ^p H^{2g-1}(A_{\CC },\QQ ),
  \wedge ^pH^1(X_{\CC },\QQ ))\; .
  $$
By Poincar\'e duality, the target of $\mu $ is the space of operators
  $$
  \Hom _{\QQ }(H^{2g-p}(A_{\CC },\QQ ),H^p(X_{\CC },\QQ ))\; .
  $$

Consider the map induced by the product in cohomology of $X$,
  $$
  \rho :\wedge ^pH^1(X_{\CC },\QQ )\lra H^p(X_{\CC },\QQ )\; ,
  $$
and the map
  $$
  \delta : H^{2g-p}(A_{\CC },\QQ )\lra \wedge^pH^{2g-1}(A_{\CC },\QQ )
  $$
which is dual to the map
  $$
  \wedge^p H^1(A_{\CC },\QQ )\lra H^p(A_{\CC },\QQ )\; ,
  $$
induced by the product in the cohomology groups of $A$. It follows from linear algebra that the $p$-th power map $\mu \circ \lambda $ sends an operator
  $$
  \psi : H^{2g-1}(A_{\CC },\QQ )\lra H^1(X_{\CC },\QQ )
  $$
to the operator
  $$
  \phi : H^{2g-p}(A_{\CC },\QQ )\lra H^p(X_{\CC },\QQ )\; ,
  $$
defined as the composition
  $$
  \phi = \rho \circ \wedge^p\psi \circ \delta \; .
  $$

In particular, if $\psi $ is the isomorphism given by the action $(\beta _{\CC })_{*,2g-1}$, the correspondence $\alpha $ acts on the cohomology groups of $A$ by the formula
  $$
  (\alpha _{\CC })_{*,2g-p}=
  \rho \circ \wedge^p((\beta _{\CC })_{*,2g-1})\circ
  \delta: H^{2g-p}(A_{\CC },\QQ )\to H^p(X_{\CC },\QQ )\; .
  $$

\bigskip

Now we are ready to construct an element $\gamma $ in the group $H^{g,g-p}(A_{\CC })$, such that $(\alpha _{\CC })_*(\gamma)\ne 0$.

\medskip

Recall that, by assumption, we have $p$ forms $\omega _1,\ldots,\omega _p$ on $X_{\CC }$ with
  $$
  \omega _1\wedge \ldots \wedge \omega _p\ne 0\; .
  $$
For each index $i$ consider the element
  $$
  \xi _i=((\beta \,_{\CC })_{*,2g-1})^{-1}(\omega_i)
  $$
in the group $H^{2g-1}(A_{\CC },\QQ )$. Since $A$ is an abelian variety, the multiplication map
  $$
  \wedge^p H^1(A_{\CC },\QQ )\lra H^p(A_{\CC },\QQ )
  $$
is an isomorphism. Hence, the dual map $\delta $ is also an isomorphism.

Take then the element
  $$
  \gamma =\delta ^{-1}(\xi _1\wedge \ldots \wedge \xi_p)
  $$
in $H^{2g-p}(A_{\CC },\QQ )$. By the above formula for
$(\alpha _{\CC })_{*,2g-p}$, we have that
  $$
  \begin{array}{rcl}
  (\alpha\,_{\CC })_*(\gamma)
  &=&
  (\rho\circ \wedge^p((\beta\,_{\CC })_{*,2g-1})\circ \delta)(\gamma) \\
  &=&
  (\rho\circ \wedge^p((\beta\,_{\CC })_{*,2g-1}))(\xi_1\wedge \ldots \wedge \xi _p) \\
  &=&
  \rho(\omega_1\wedge\ldots\wedge\omega_p) \\
  &=&
  \omega_1\wedge\ldots\wedge\omega_p \\
  &\neq &
  0 \\
  \end{array}
  $$
The differential form $\omega _1\wedge \ldots \wedge \omega _p$ belongs to the group $H^{p,0}(X_{\CC })$. Since the correspondence $\alpha $ is of degree $p-g$, the map $(\alpha \,_{\CC })_*$ shifts the indices in the Hodge filtration by $p-g$. Therefore, $\gamma $ sits in $H^{g,g-p}(A_{\CC })$.
Then, by Theorem~\ref{non-triviality}, the cycle class $\alpha _{\eta }$ is non-trivial.
\end{pf}

\medskip

\begin{remark}
{\rm Applying a similar argument like in Remark \ref{lower}, one can also show the following. Let $X$ be an irreducible smooth projective variety over $k$. As above, suppose there exist $p$ differential forms
  $$
  \omega _1,\dots,\omega _p\in
  H^0(X_{\CC },\Omega^1_{X_{\CC }})=H^{1,0}(X_{\CC })\; ,
  $$
such that
  $$
  \omega _1\wedge \dots \wedge \omega _p\ne 0\; .
  $$
Take a smooth linear section $T$ of codimension $g-p$ in $A$, let $f:T\hra A$ be the corresponding closed embedding, and let
  $$
  \beta =
  (f\times \id _X)^*(\alpha )\in CH^p(T\times X)=
  CH^0(T,X)\; ,
  $$
be the restriction of the class $\alpha $ on $T\times X$. Let also $\xi $ be the generic point of $T$. Then $\beta _{\xi }$ is a non-zero element in the Abel-Jacobi kernel $T^p(X_{k(T)})$.}
\end{remark}

\bigskip

Finally, we would like to take a closer look at the condition in Theorem \ref{oneforms}. Let $\omega $ be a holomorphic $1$-form on $X$, i.e. a section of the cotangent bundle $T_X^*$ on the complex manifold $X_{\CC }$. For a closed point $x\in X_{\CC }$ denote by $\omega|_x$ the value of the differential form $\omega $ at $x$. Thus, $\omega|_x$ is an element in the cotangent space $T^*_{X,x}$ to $X_{\CC }$ at $x$. The condition $\omega _1\wedge \ldots \wedge \omega _p\neq 0$ in Theorem \ref{oneforms} is equivalent to say that there exists a closed point $x\in X_{\CC }$, such that
  $$
  \omega _1|_x\wedge \ldots \wedge \omega _p|_x\neq 0\; .
  $$
It follows from linear algebra that $\omega _1|_x\wedge \ldots \wedge \omega _p|_x\neq 0$ if and only if the vectors
  $$
  \omega _1|_x,\ldots ,\omega _p|_x
  $$
are linearly independent in the cotangent space $T^*_{X,x}$. The next lemma shows that latter thing is equivalent to the condition that the image of the Albanese mapping $X\to \Alb(X)$ is at least $p$-dimensional\footnote{This lemma seems to be well known for experts, but we prove it here for the convenience of the reader}.

\begin{lemma}
\label{lemma-exampforms}
Let $X$ be an irreducible smooth projective variety $X$ over $k$. Then it has $p$ holomorphic one-forms $\omega _1,\ldots ,\omega _p$ in $H^0(X_{\CC },\Omega ^1_{X_{\CC }})=H^{1,0}(X_{\CC })$, such that
  $$
  \omega _1\wedge \ldots \wedge \omega _p\neq 0\; ,
  $$
if and only if the dimension of the image of the Albanese map $X\to \Alb (X)$ is greater or equal than $p$.
\end{lemma}

\begin{pf}
For short, let $B=\Alb (X)$, and let
  $$
  f:X\lra B
  $$
be the Albanese mapping.

First assume that there are $p$ forms $\omega_1,\ldots,\omega_p$ as above. Denote by $Y$ the image of the Albanese map $f:X\to B$. Since $X$ is irreducible, $Y$ is irreducible too. Let $U$ be a non-empty open subset in $Y$, such that $U$ is a smooth variety over $k$, and put
  $$
  V=f^{-1}(U)\; .
  $$
Notice that $V$ is a non-empty open subset in $X$, because the natural map $X\to Y$ is surjective. Let
  $$
  i:V\lra X
  $$
be the corresponding open embedding. Denote by
  $$
  g:V\lra U
  $$
the restriction of the natural map $X\to Y$ to $V$, and let
  $$
  j:U\lra B
  $$
be the composition $U\to Y\to B$.

Let $W$ be the set of closed points $x\in X$ such that
   $$
   \omega_1|_x\wedge\ldots\wedge\omega_p|_x\ne 0\; .
   $$
Then $W$ is an open subset in $X$ and, moreover, $W$ is non-empty because $\omega_1\wedge\ldots\wedge\omega_p\neq 0$.

Since $X$ is irreducible, the intersection of two non-empty open subsets $W\cap V$ is non-empty. Thus, there exists a closed point $x\in V$ with $\omega_1|_x\wedge\ldots\wedge\omega_p|_x\neq 0$. Therefore, the $p$-form $i_{\CC }^*(\omega_1\wedge\ldots\wedge \omega_p)$ on $V$ is non-zero too.

As all holomorphic $1$-forms on $X_{\CC }$ are obtained by pulling-back via $f_{\CC }^*$ of holomorphic $1$-forms on $B_{\CC }$, there are $p$ holomorphic $1$-forms $\xi_1,\ldots,\xi_p$ on $B_{\CC }$, such that
  $$
  f^*_{\CC }\,\xi_i=\omega_i
  $$
for all $i$.

Since $f\circ i=j\circ g$, one has
  $$
  g^*_{\CC }j^*_{\CC }(\xi_1\wedge \ldots \wedge \xi _p)=
  i^*_{\CC }f^*_{\CC }(\xi_1\wedge\ldots\wedge \xi_p)\; .
  $$
As
  $$
  f^*_{\CC }(\xi_1\wedge\ldots\wedge\xi_p)=\omega_1\wedge\ldots\wedge\omega_p
  $$
and the form $i^*_{\CC }(\omega_1\wedge\ldots\wedge\omega_p)$
does not vanish,
  $$
  g^*_{\CC }j^*_{\CC }(\xi_1\wedge\ldots\wedge\xi_p)\neq 0\; .
  $$
This implies that
  $$
  j^*_{\CC }(\xi_1\wedge\ldots\wedge\xi_p)\neq 0\; .
  $$

Thus, there exists a non-zero $p$-form, namely $j^*_{\CC }(\xi_1\wedge\ldots\wedge\xi_p)$, on the smooth variety $U$. Hence, $\dim (U)\geq p$. As $U$ is dense in $Y$, we have that $\dim(Y)=\dim(U)$, whence $\dim (Y)\geq p$.

\medskip

Suppose now that $\dim (f(X))=\dim (Y)\geq p$. Take a smooth closed point $x$ in $Y$. Choose any codimension $p$ linear subspace
  $$
  E\subset T_xB
  $$
that intersects with $T_xY$ only by zero. Then there are $p$ linearly independent vectors
  $$
  l_1,\ldots,l_p
  $$
in the cotangent space $T_x^*B$, such that
  $$
  l_i|_E=0\; ,
  $$
where elements in the cotangent space $T_x^*B$ are considered as linear functionals on the tangent space $T_xB$.

Since $B$ is an abelian variety, the tangent bundle on $B$ is trivial. Therefore, there are $p$ forms
  $$
  \xi _1,\ldots ,\xi _p
  $$
in $H^0(B,\Omega ^1_B)$, such that
  $$
  (\xi _1)|_x=l_1,\ldots ,(\xi _p)|_x=l_p\; .
  $$
We now have that
  $$
  (\xi _1)|_x\wedge \ldots \wedge (\xi _p)|_x\neq 0\; .
  $$
It follows from linear algebra that the skew covector
  $$
  l_1\wedge \ldots \wedge l_p
  $$
remains non-zero when being restricted to the subspace $T_xY$ in $T_xB$. Therefore, the restriction of the form $\xi _1\wedge \ldots \wedge \xi _p$ from $B$ to $Y$ is non-zero at the point $x$. Thus, the restriction of the whole form $\xi_1\wedge\ldots\wedge\xi_p$ from $B$ to $Y$ is non-zero.

Let
  $$
  \omega _i=f^*\xi _i
  $$
for each index $i$. Since the morphism $X\to Y$ is surjective, the form
  $$
  \omega _1\wedge \ldots \wedge \omega _p
  $$
is non-zero too.
\end{pf}

\medskip

Now one can make our construction to be more concrete. Namely, take an abelian variety $B$ over $k$. Then take any subvariety $Y$ in $B$, such that $\dim (Y)\geq p$. Consider a smooth projective variety $X$ that admits a surjective morphism onto $Y$. Then the dimension of the image of the Albanese map $f:X\to \Alb (X)$ is greater or equal than $p$. Applying Theorem \ref{oneforms}, we get a non-zero element in the Abel-Jacobi kernel $T^p(X)$ of the variety $X$.

\bigskip

\begin{small}

\end{small}

\bigskip

\bigskip

\begin{small}

{\sc Steklov Mathematical Institute, Gubkina str. 8,
119991, Moscow, Russia}

\end{small}

\medskip

\begin{footnotesize}

{\it E-mail address}: {\tt gorchins@mi.ras.ru}

\end{footnotesize}

\bigskip

\begin{small}

{\sc Department of Mathematical Sciences, University of Liverpool,
Peach Street, Liverpool L69 7ZL, England, UK}

\end{small}

\medskip

\begin{footnotesize}

{\it E-mail address}: {\tt vladimir.guletskii@liverpool.ac.uk}

\end{footnotesize}

\end{document}